\documentclass[openany, a4paper, 11pt]{article}

\usepackage{amsfonts}

 \usepackage{mathrsfs,amsfonts,amsmath}
 \usepackage{color}

 \makeatletter
 \@addtoreset{equation}{section}
 \makeatother
 \newtheorem{theorem}{Theorem}[section]
 \newtheorem{definition}{Definition}[section]
 \newtheorem{hypothesis}{Hypothesis}[section]
 \newtheorem{lemma}{Lemma}[section]
 \newtheorem{proposition}{Proposition}[section]
 \newtheorem{corollary}{Corollary}[section]
 \newtheorem{remark}{Remark}[section]
 \newtheorem{example}{Example}[section]

 \def\beqlb{\begin{eqnarray}}\def\eeqlb{\end{eqnarray}}
 \def\beqnn{\begin{eqnarray*}}\def\eeqnn{\end{eqnarray*}}

 \def\ar{\!\!\!&}

 \def\mbb{\mathbb}

 \def\qed{\hfill$\Box$\medskip}

\newcommand{\bcen}{\begin{center}}
\newcommand{\ecen}{\end{center}}
\newcommand{\bgeqn}{\begin{equation}}
\newcommand{\edeqn}{\end{equation}}

\def\dz{\delta}
\def\lz{\lambda}
\def\oz{\omega}

\def\wz{\wedge}

\def\ra{\rangle}
\def\la{\langle}

\begin{document}

\vskip1.0cm

\centerline{\large\textbf{Discontinuous Superprocesses with
Dependent Spatial Motion }\footnote{ Supported by NSFC (No.
10121101)}}

\bigskip

\centerline{ Hui He\footnote{ \textit{E-mail address:} {
h$_{-}$hui$_{-}$math@mail.bnu.edu.cn }} }

\medskip

\centerline{Laboratory of Mathematics and Complex Systems, }

\centerline{ School of Mathematical Sciences, Beijing Normal
University,}

\smallskip

\centerline{ Beijing 100875, People's Republic of China}

\bigskip

{\narrower{\narrower{\narrower{\narrower

\begin{abstract}
We construct a class of discontinuous superprocesses with dependent
spatial motion and general branching mechanism. The process arises
as the weak limit of critical interacting-branching particle systems
where the spatial motions of the particles are not independent. The
main work is to solve the martingale problem. When we turn to the
uniqueness of the process,  we generalize the localization method
introduced by [D.W. Stroock, Diffusion processes associated with
L\'{e}vy generators, Z. Wahrscheinlichkeitstheorie und Verw.
Gebiete, 32(1975) 209--244] to the measure-valued context. As for
existence, we use particle system approximation and a perturbation
method. This work generalizes the model introduced in [D.A. Dawson,
Z. Li, H. Wang, Superprocesses with dependent spatial motion and
general branching densities, Electron. J. Probab. 6(2001), no.25, 33
pp. (electronic)] where quadratic branching mechanism was
considered. We also investigate some properties of the process.
\end{abstract}

\smallskip


\smallskip

\smallskip

\noindent\textbf{Abbreviated Title:} Discontinuous superprocesses
\par}\par}\par}\par}

\bigskip\bigskip

\section{Introduction}
Notation: For reader's convenience, we introduce here our main
notation. Let $\hat{\mathbb{R}}$ denote the  one-point
compactification of $\mathbb{R}$. Let $\hat{\mbb{R}}^n$ denote the
$n$-fold Cartesian product of $\hat{\mbb{R}}$. Let $M(\mbb{R})$
denote the space of finite measure endowed with topological of weak
convergence. We denote by $\lz^n$ the Lebesgue measure on
$\mathbb{R}^n$. Given a topological space $E$, let $\mathfrak{B}(E)$
denote borel $\sigma$-algebra on $E$. Let $B(E)$ denote the set of
bounded measurable functions on $E$ and
  let $C(E)$ denote its subset comprising of bounded continuous
  functions. Let $\hat{C}(\mathbb{R}^n)$ be the space of continuous
   functions  on $\mbb R^n$ which vanish at infinity and
  let $C_c^{\infty}(\mathbb{R}^n)$ be functions with compact support
 and bounded continuous derivatives of any order.
  Let $C^2(\mathbb{R}^n)$ denote the set of functions in
$C(\mathbb{R}^n)$ which is twice continuously differential functions
with bounded derivatives up to the second order. Let $C_c^2(\mbb
R^n)$ denote the set of functions in $C^2(\mbb R^n)$ with compact
support. Let $\hat{C}^2(\mathbb{R}^n)$ be the subset of
$C^2(\mathbb{R}^n)$ of functions that together with their
derivatives up to the second order vanish  at infinity. \\
\noindent
Let
$$
C_{\partial}^2(\mbb{R}^n)=\{f+c: c\in\mbb{R} \textrm{ and } f\in
\hat{C}^2({\mbb{R}^n})\}
$$
and
$$C_0^2(\mbb{R}^n)=\{f: f\in{C}^2_{\partial}(\mbb{R}^n)
\textrm{ and }(1+|x|^2)D^{\alpha}f(x)\in
\hat{C}(\mbb{R}^n),~\alpha=1,2\},
$$
where $D^1f=\sum\limits_{i=1}^n|\partial f/\partial x_i|$ and
$D^2f=\sum\limits_{i,j=1}^{n}|\partial^2f/\partial x_i\partial x_j
|$. We use the superscript ``+'' to denote the subsets of
non-negative elements of the function spaces, and ``++'' is used to
denote the subsets of non-negative elements bounded away from zero,
e.g., $B(\mbb R^n)^+$, $C(\mbb R^n)^{++}$. Let $f^i$ denote the
first order partial differential derivatives of the function
$f(x_1,\cdots,x_n)$ with respect to $x_i$ and
 let  $f^{ij}$ denote the second
order partial differential derivatives of the function
$f(x_1,\cdots,x_n)$ with respect to $x_i$ and $x_j$.
 We denote by $C\left([0,\infty),E\right)$ the space of continuous paths taking values in
 $E$.
 Let $D\left([0,\infty),E\right)$ denote the Skorokhod space of c\`{a}dl\`{a}g paths taking values
 in $E$. For $f\in C(\mathbb{R})$ and $\mu\in M(\mathbb{R})$ we shall write
$\langle f,\mu\rangle$ for $\int fd\mu$.

\bigskip

A class of  \textit{superprocesses with dependent spatial motion}
(SDSM) over the real line $\mathbb{R}$ were introduced and
constructed in \cite{[W97],[W98]}. A generalization of the model was
then given in \cite{[DLW01]}. We first briefly describe the model
constructed in \cite{[DLW01]}. Suppose that $c\in C^2(\mathbb{R})$
 and $h\in C(\mathbb{R})$ is square-integrable. Let
 \bgeqn \label{f1.1}
\rho(x)=\int_{\mathbb{R}}h(y-x)h(y)dy,
 \edeqn
and $a(x)=c(x)^2+\rho(0)$ for $x\in \mathbb{R}.$ We assume in
addition that  $\rho\in C^2(\mbb R)$ and $|c|$ is bounded away from
zero. Let $\sigma$
 be a nonnegative function in $C^2(\mathbb{R})$ and can be extended
continuously to $\hat{\mathbb{R}}$. Given a finite measure $\mu$ on
$\mathbb{R}$, the SDSM with parameters $(a, \rho, \sigma)$ and
initial state $\mu$ is the unique solution of the
$(\mathcal{L},\mu)$-martingale problem, where
 \bgeqn
 \label{f1.2}
 \mathcal{L}F(\mu):=\mathcal{A}F(\mu)+\mathcal{B}F(\mu),
 \edeqn
 \begin{eqnarray}
\label{f1.3}
\mathcal{A}F(\mu)\ar:=\ar\frac{1}{2}\int_{\mathbb{R}}a(x)\frac{d^2}{dx^2}
                    \frac{\dz F(\mu)}{\dz\mu(x)}\mu(dx)\cr&&
+\frac{1}{2}\int_{\mathbb{R}^2}\rho(x-y)\frac{d^2}{dxdy}
                  \frac{\dz^2F(\mu)}{\dz\mu(x)\dz\mu(y)}\mu(dx)\mu(dy),
\end{eqnarray}
 \bgeqn \label{f1.4}
\mathcal{B}F(\mu):=\frac{1}{2}\int_\mathbb{R}\sigma(x)
\frac{\dz^2F(\mu)}{\dz\mu(x)^2}\mu(dx),
 \edeqn
for some bounded continuous functions $F(\mu)$ on $M(\mathbb{R})$.
The variational derivative is defined by
\begin{equation}
\label{f1.5} \frac{\dz
F(\mu)}{\dz\mu(x)}=\lim_{r\rightarrow0^+}\frac{1}{r}[F(\mu+r\dz_x)-F(\mu)],\textrm{\
\ \ } x\in \mathbb{R},
\end{equation}
if the limit exists and  $\dz^2F(\mu) / \dz\mu(x)\dz\mu(y)$ is
defined in the same way with $F$ replaced by $(\dz F/\dz\mu(y))$ on
the right hand side. Clearly, the SDSM reduces to a usual critical
 Dawson-Watanabe superprocess if $h(\cdot)\equiv0$ (see \cite{[D93]}).
 A general SDSM arises as the weak limit of critical
interacting-branching particle systems. In contrast to the usual
branching particle system, the spatial motions of the particles in
the interacting-branching particle system are \textit{not}
independent. The spatial motions of the particles can be described
as follows. Suppose that $\{ W(t,x): x\in \mathbb{R}, t\geq0\}$ is
space-time white noise based on Lebesgue measure, the common noise,
and $\{B_i(t):t\geq0,i=1,2,\cdots\}$ is a family of independent
standard Brownian motions, the individual noises, which are
independent of $\{W(t,x): x\in\mathbb{R}\}$. The migration of a
particle in the approximating system with label $i$ is defined by
the stochastic equations
 \bgeqn
  \label{1.6}
dx_i(t)=c(x_i(t))dB_i(t)+\int_\mathbb{R}h(y-x_i(t))W(dt,dy),\textrm{\
\ \ }t\geq0,~~i=1,2,\cdots,
 \edeqn
where $W(dt,dy)$ denotes the time-space stochastic integral relative
to $\{W_t(B)\}$. For each integer $m\geq1$,
$\{(x_1(t),\cdots,x_m(t)):t\geq0\}$ is an $m$-dimensional diffusion
process which is generated by the differential operator
 \bgeqn
 G^m:=\frac{1}{2}\sum_{i=1}^ma(x_i)\frac{\partial^2}{\partial
x_i^2}+ \frac{1}{2}\sum_{i,j=1,i\neq
j}^m\rho(x_i-x_j)\frac{\partial^2}{\partial x_i\partial x_j}.
 \edeqn
In particular, $\{x_i(t):t\geq0\}$ is a one-dimensional diffusion
process with generator $G:=(a(x)/2)\Delta$. Because of the
exchangeability, a diffusion process generated by $G^m$ can be
regarded as an interacting particle system or a measure-valued
process. Heuristically, $a(\cdot)$ represents the speed of the
particles and $\rho(\cdot)$ describes the interaction between them.
The diffusion process generated by $\mathcal{A}$ arises as the high
density limit of a sequence of interacting particle systems
described by $(\ref{1.6})$; see Wang \cite{[W97],[W98]} and Dawson
\emph{et al} \cite{[DLW01]}.  There are at least two different ways
to look at the SDSM. One is as a superprocess in random environment
and the other as an extension of models of the motion of the mass by
stochastic flows (see \cite{[MX01]}). Some other related models were
introduced and studied in Skoulakis and Adler \cite{[SA01]}. The
SDSM possesses properties very different from those of the usual
Dawson-Watanabe superprocess. For example, a Dawson-Watanabe
superprocess in $M(\mathbb{R})$ is usually absolutely continuous
whereas the SDSM with $c(\cdot)\equiv0$ is purely atomic; see Konno
and Shiga \cite{[KS88]} and  \cite{[DL03], [W02]}, respectively.

To best of our knowledge, in all of the work which considered the
SDSM and related models only continuous processes have been
introduced and studied. In this paper, we construct a class of
discontinuous superprocesses with dependent spatial motion.  A
modification of the above martingale problem is to replace operator
$\mathcal{B}$ in ($\ref{f1.2}$) by
 \begin{eqnarray}
 \label{f1.8}
\mathcal{B}F(\mu)\ar=\ar\frac{1}{2}\int_\mathbb{R}\sigma(x)
\frac{\dz^2F(\mu)}{\dz\mu(x)^2}\mu(dx)\cr
&&+\int_\mathbb{R}\mu(dx)\int_0^\infty
\left(F(\mu+\xi\dz_x)-F(\mu)-\frac{\dz
F(\mu)}{\dz\mu(x)}\xi\right)\gamma(x,d\xi),
 \end{eqnarray}
whose coefficients satisfy:
\begin{enumerate}
\item[(i)] $\sigma\in C_{\partial}^2(\mbb{R})^+$,
\item[(ii)] $\gamma(x,d\xi)$ is a kernel from $\mathbb{R}$
to $(0, +\infty)$ such that
$\sup\limits_x[\int_0^{+\infty}\xi\wedge\xi^2\gamma(x,d\xi)]<+\infty$,
\item[(iii)] $\int_{\Gamma} \xi\wedge\xi^2\gamma(x,d\xi)
\in C^2_{\partial}(\mbb{R})$ for each $\Gamma\in
\mathfrak{B}((0,\infty))$.
\end{enumerate}

\noindent A Markov process generated by $\mathcal{L}$ is a
measure-valued  branching process with branching mechanism  given by
$$
 \Psi(x,z):=\frac{1}{2}\sigma(x)z^2
  +\int_0^{\infty}(e^{-z\xi}-1+z\xi)\gamma(x,d\xi).
$$
 This process is naturally called a
\textsl{superprocess with dependent spatial motion (SDSM) }with
parameters $( a,\rho,\Psi)$. This modification is related to the
recent work of Dawson \emph{et al} \cite{[DLW01]}, where it was
assumed that $\gamma(x,d\xi)=0$.  Though our model is an extension
of the model introduced in Wang \cite{[W97],[W98]} and Dawson
\emph{et al} \cite{[DLW01]}, the construction of our model  differ
from theirs. We describe our approach to the construction of our
model in the following.

The main work of this paper is to solve the
$(\mathcal{L},\mu)$-martingale problem. As for uniqueness, following
the idea of Stroock \cite{[S75]} a localization procedure is
developed. Therefore, we do not consider the $({\cal
L},\mu)$-martingale problem directly. Instead,  we will first solve
the $({\cal L'},\mu)$-martingale problem, where
\begin{eqnarray}
\label{f1.9}
\mathcal{L}'F(\mu):=\mathcal{A}F(\mu)+\mathcal{B}'F(\mu),
\end{eqnarray}
\begin{eqnarray}
 \mathcal{B}'\ar:=\ar\frac{1}{2}\int_\mathbb{R}\sigma(x)
\frac{\dz^2F(\mu)}{\dz\mu(x)^2}\mu(dx)-\int_{\mathbb{R}}\mu(dx)\int_l^{\infty}\frac{\dz
F(\mu)}{\dz\mu(x)}\xi\gamma(x,d\xi)\cr
&&+\int_\mathbb{R}\mu(dx)\int_0^l
\left(F(\mu+\xi\dz_x)-F(\mu)-\frac{\dz
F(\mu)}{\dz\mu(x)}\xi\right)\gamma(x,d\xi).
\end{eqnarray}
and we make the convention that
 \begin{eqnarray*}
\int_0^l = \int_{(0,l)}
 \qquad\mbox{and}\qquad
\int_l^\infty = \int_{[l,\infty)}
 \end{eqnarray*}
for $0<l <\infty$. We regard the $({\cal L'},\mu)$-martingale
problem as the `killed' martingale problem. We shall see that the
Markov process associated with the `killed' martingale problem also
arises as high density limit of a sequence of interacting-branching
particle system and it is an SDSM with branching mechanism given by
$$
 \Psi_0(x,z):=\frac{1}{2}\sigma(x)z^2+\int_l^{\infty}\xi\gamma(x,d\xi)z
  +\int_0^l(e^{-z\xi}-1+z\xi)\gamma(x,d\xi).
$$
 It is easy to see from the branching mechanism
that  the process is a subcritical branching process with all `big'
jumps such that the jump size is larger than $l$ been `killed'.
 We will use duality method to show
the uniqueness of the `killed' martingale problem. We shall
construct a dual process and show its connection with  the solutions
of the `killed' martingale problem which gives the uniqueness. When
we establish the dual relationship, we point out that there exists a
gap in the proof of establishing the dual relationship in
\cite{[DLW01]}; see Remark $\ref{gap}$ in Section 2 of this paper
for details. Then  a localization argument is developed  to show
that if the $({\cal L}', \mu)$ martingale problem is well-posed then
uniqueness holds for the $({\cal L},\mu)$-martingale problem. The
argument  consists of three parts.

In the first part, we show that each solution of the $({\cal
L},\mu)$-martingale problem , say $X$, behaves the same as the
solution of the killed martingale problem until it has a `big jump'
whose jump size is larger than $l$. Intuitively, one can think of
the branching particle system as follows. In the branching particle
system corresponding to  the $({\cal L},\mu)$-martingale problem, if
a particle dies and it leaves behind a large number of offsprings,
say more than 500, which always be regarded as a `big jump' event,
we kill all its offsprings. Then we get a new branching particle
system and before the jump event happens the two systems are the
same. The evolution of the new particle system represents the
behavior of the solution to the `killed' martingale problem.  It is
clear that if the original branching particle system is a critical
system, the new particle system is a subcritical branching system.
Since the `killed' martingale problem is well-posed, $X$ is uniquely
determined before it has a `big jump'. Next, we show that when a
`big jump' event happens, the jump size is uniquely determined. This
conclusion is not surprising either. Given a branching mechanism, in
a branching particle system, when a particle dies, the distribution
of its offspring number is uniquely determined by the position of
the particle itself (we assume that the branching mechanism is
independent of time). Thus we can find a predictable representation
for the jump size. According to the argument in the first part, we
see the jump size is uniquely determined. At last, we can prove by
induction that the distribution of $X$ is uniquely determined, since
after the first `big jump' event happens, $X$ also behaves the same
as the solution of the `killed' martingale problem until the second
`big jump' event happens. Before we use the localization procedure,
we follow  an argument taken from El-Karoui and Roelly-Coppoletta
\cite{[ER91]} to decompose each solution of the $({\cal
L},\mu)$-martingale problem into a continuous part and a purely
discontinuous part. We will use this argument again when we show the
existence of solutions to the $({\cal L},\mu)$-martingale problem;
see next two paragraphs.

When we turn to the existence  we also first consider the existence
of the `killed' martingale problem. Although the solution of the
`killed' martingale problem is also an SDSM which arises as high
density limit of a sequence of interacting-branching particle
systems,  in order to deduce the martingale formula the techniques
developed in Wang \cite{[W97],[W98]} and Dawson \emph{et al}
\cite{[DLW01]} can not be used directly because of the third item in
the branching mechanism $\Psi_0$. We will use the martingale
decomposition and special semi-martingale's representation to get
the desired result. Our approach is stimulated by  El-Karoui and
Roelly-Coppoletta \cite{[ER91]}, who considered the martingale
problem of the usual Dawson-Watanabe superprocess. We briefly
describe the main idea in next paragraph.

First, a sequence of subcritical branching particle systems is
constructed. Let $X=(X_t)_{t\geq0}$ denote a limit of the particle
systems. Then we derive the special semi-martingale property of
 $\{\exp\{-\la\phi,X_t\ra\}:t\geq0\}$ with $\phi$ bounded away
 from zero by using particle system approximation,
 and obtain a representation for this semi-martingale. This
 approach is \emph{different} from that of \cite{[ER91]}, where log-laplace
 equation was used to deduce the semi-martingale property.
 Next, we consider
 an integer-valued random measure $N(ds,d\nu)=\sum_{s>0}1_{\{\Delta
 X_s\neq0\}}\delta_{(s,\Delta X_s)}(ds,d\nu)$ and by an approximation procedure
 we can show
 \bgeqn
 \label{f1.11}
 M_t(\phi):=\la\phi,X_t\ra-\la\phi, X_0\ra-\frac{1}{2}\int_0^t\la
 a\phi'',X_s\ra ds+\int_0^tds\la\int_l^{\infty}\xi\gamma(\cdot,d\xi)\phi,X_s\ra
 \edeqn
 is square-integrable martingale which can be decompose
 into a continuous martingale $\{M_t^c(\phi):t\geq0\}$ and a purely discontinuous
 martingale $\{M_t^d(\phi):t\geq0\}$.  We have
 \bgeqn
 \label{f1.12}
 \la \phi,X_t\ra=\la\phi, X_0\ra+\frac{1}{2}\int_0^t\la
 a\phi'',X_s\ra ds+ M_t^c(\phi)+M_t^d(\phi)-
 \int_0^tds\la\int_l^{\infty}\xi\gamma(\cdot,d\xi)\phi,X_s\ra,
 \edeqn
 and  $M^d(\phi)$  can be represented as a stochastic
 integral with respect to the corresponding martingale measure of $N(ds,d\nu)$.
 This argument is also different from the argument of
 \cite{[ER91]}, where according to the semi-martingale property of
 $\{\exp\{-\la\phi,X_t\ra\}:t\geq0\}$ only semi-martingale property of
 $\{\la\phi,X_t\ra:t\geq0\}$ with $\phi$ bounded away from zero
 was derived.
  By the martingale decomposition ($\ref{f1.12}$) we can obtain
 another representation
 for semi-martingale
 $\{\exp\{-\la\phi,X_t\ra\}:t\geq0\}$.
 By identifying  two  representations for
 $\{\exp\{-\la\phi,X_t\ra\}:t\geq0\}$ mentioned above,
 we know the explicit form of the quadratic variation process
 of  $\{M^c_t(\phi):t\geq0\}$ and the
 compensator of the random measure $N(ds,d\nu)$.
 Then we can deduce $X$ satisfies the martingale formula for  the
 $({\cal L}',\mu)$-martingale problem.
  At last by a  perturbation method we show the existence of the $({\cal
 L},\mu)$-martingale problem.

 The remainder of the paper is organized as follows. In Section 2,
 we  first introduce the `killed' martingale problem
  and define a  dual process and investigate its
 connection to the solutions of the `killed' martingale problem which gives
 the uniqueness of the `killed' martingale problem. Then we deduce that the
 uniqueness holds for the $({\cal L},\mu)$-martingale problem.
 In Section 3, we first give a formulation
 of the system of branching particles with dependent spatial motion
and obtain the existence of the solution of the `killed' martingale
problem by taking high density limit of particle systems. Then a
perturbation argument is used to show the existence of the $({\cal
L},\mu)$-martingale problem. We compute the first and second order
moment formulas of the process in Section 4.

\begin{remark}
\label{rem1.1}By Theorem 8.2.5 of \cite{[EK86]},  the closure of
$\{(f,G^mf):f\in C_c^{\infty}(\mathbb{R}^m)\}$ which we still denote
by $G^m$ is single-valued and generates a Feller semigroup
$(P_t^m)_{t\geq0}$ on $\hat{C}(\mathbb{R}^m)$. Note that this
semigroup is given by a transition function and can therefore be
extended to all of $B(\mathbb{R}^m)$. We also have that $(1,0)$ is
in the bp-closure of $G^m$.
\end{remark}

\section{Uniqueness}
\subsection{Killed martingale problem}
In this section, we first introduce  the \textsl{killed martingale
problem} for the SDSM and show the uniqueness holds for the killed
martingale problem.

\noindent \begin{definition}Let $\mathcal {D}(\mathcal
{L})=\bigcup_{m=0}^{\infty}\left\{F(\mu)
=f(\la\phi_1,\mu\ra,\cdots,\la\phi_m,\mu\ra),~
f\in C_0^2(\mathbb{R}^m),~\{\phi_i\}\subset
C_c^2(\mathbb{R})^+\right\}$. For $\mu\in M(\mathbb{R})$ and  an
$M(\mathbb{R})$-valued c\`{a}dl\`{a}g process $\{X_t:t\geq0\}$, we
say $X$ is a solution of the $(\mathcal {L},\mu)$-{martingale
problem} if $X_0=\mu$ and
 \bgeqn\label{f2.1}
  F(X_t)-F(X_0)-\int_0^t\mathcal {L}F(X_s)ds,\textrm{\ \ \ } t\geq0,
 \edeqn
is a local martingale for each $F\in\mathcal {D}(\mathcal {L})$ and
for $l>1$, we say $X$ is a solution of the $(\mathcal
{L}',\mu)$-martingale problem if $X_0=\mu$ and \bgeqn \label{f2.2}
  F(X_t)-F(X_0)-\int_0^t\mathcal {L}'F(X_s)ds,\textrm{\ \ \ } t\geq0,
\edeqn is a local martingale for each $F\in\mathcal {D}(\mathcal
{L})$. \end{definition}
Let ${\cal D}_0({\cal
L})=\bigcup_{m=0}^{\infty}\{f(\la\phi_1,\mu\ra,\cdots,\la\phi_m,\mu\ra),~
f\in C_0^2(\mathbb{R}^m),~\{\phi_i\}\subset C^2(\mathbb{R})^{++}\}$.
Note that for $F(\mu)\in{\cal D}_0(\cal {L})\cup{\cal D}(\cal {L})$,
 \begin{eqnarray}
 \label{f2.3}
 \mathcal{A}F(\mu)&=&\frac{1}{2}\sum_{j=1}^m
 f^i(\la\phi_1,\mu\ra,\cdots,\la\phi_m,\mu\ra)\la
 a\phi''_i,\mu\ra\cr
 &&+\frac{1}{2}\sum_{i,j=1}^m
 f^{ij}(\la\phi_1,\mu\ra,\cdots,\la\phi_m,\mu\ra)
 \int_{\mathbb{R}^2}\rho(x-y)\phi_i'(x)\phi_j'(y)\mu^2(dxdy),
 \end{eqnarray}
\begin{eqnarray}
 \label{f2.4}
 \mathcal{B}F(\mu)&=&\frac{1}{2}\sum_{i,j=1}^m
 f^{ij}(\la\phi_1,\mu\ra,\cdots,\la\phi_m,\mu\ra)
 \la\sigma\phi_i\phi_j,\mu\ra\cr
 &&+\int_{\mathbb{R}}\mu(dx)\int_0^{\infty}\{
 f(\la\phi_1,\mu\ra+\xi\phi_1(x),\cdots,
 \la\phi_m,\mu\ra+\xi\phi_m(x))\cr
 &&~~~-f(\la\phi_1,\mu\ra,\cdots,\la\phi_m,\mu\ra)
  -\xi\sum_{i=1}^m
  f^i(\la\phi_1,\mu\ra,\cdots,\la\phi_m,\mu\ra)\phi_i(x)\}
  \gamma(x,d\xi)
 \end{eqnarray}
 and
 \begin{eqnarray}
 \label{f2.5}
 \mathcal{B}'F(\mu)&=&\frac{1}{2}\sum_{i,j=1}^m
 f^{ij}(\la\phi_1,\mu\ra,\cdots,\la\phi_m,\mu\ra)
 \la\sigma\phi_i\phi_j,\mu\ra\cr
 &&-\int_{\mbb R}\mu(dx)\int_l^{\infty}\xi\gamma(x,d\xi)\sum_{i=1}^m
  f^i(\la\phi_1,\mu\ra,\cdots,\la\phi_m,\mu\ra)\phi_i(x)\cr
 &&+\int_{\mathbb{R}}\mu(dx)\int_0^{l}\{
 f(\la\phi_1,\mu\ra+\xi\phi_1(x),\cdots,
 \la\phi_m,\mu\ra+\xi\phi_m(x))\cr
 &&~~~-f(\la\phi_1,\mu\ra,\cdots,\la\phi_m,\mu\ra)
  -\xi\sum_{i=1}^m
  f^i(\la\phi_1,\mu\ra,\cdots,\la\phi_m,\mu\ra)\phi_i(x)\}
  \gamma(x,d\xi).
 \end{eqnarray}
Thus for every $F\in {\cal {D}}_0({\cal {L}})$, both ${\cal L}F$ and
${\cal L'}F$ are bounded functions on $M(\mbb R)$.
\begin{remark}
\label{rem2.1} Let $h\in C_c^2(\mbb R^m)$ satisfy $1_{B(0,1)}\leq
h\leq 1_{B(0,2)}$ and $h_k(x)=h(x/k)\in C_c^2(\mbb R^m)$. Then for
each $\phi\in C^2(\mbb R)^{++}$, it can be approximated by $\{\phi
h_k\}\subset C_c^2(\mbb R)^{+}$ in such a way that not only $\phi$
but its derivatives up to second order are approximated boundedly
and pointwise. Therefore when $X$ is a solution of $(\cal
L,\mu)$-martingale problem (or $(\cal L',\mu)$-martingale problem),
(\ref{f2.1}) (or (\ref{f2.2})) is  a martingale for $F\in {\cal
D}_0(\cal L)$. On the other hand, for every $\phi\in C_c^2(\mbb
R)^+$, we can approximate $\phi$ by $\{\phi+1/n\}\subset
C^2_{\partial}(\mbb R)^{++}\subset C^2(\mbb R)^{++}$ in the same
way. Thus if (\ref{f2.1}) (or (\ref{f2.2}))  is a martingale for
every $F\in {\cal D}_0(\cal L)$, it is a local martingale for every
$F\in {\cal D}(\cal L)$. We shall see that any solution of the
$({\cal L}',\mu)$-martingale problem has bounded moment of any
order. Thus if $X$ is a solution of the  $(\cal L',\mu)$-martingale
problem, (\ref{f2.2}) is  a martingale for every $F\in{\cal
D}_0(\cal {L})\cup{\cal D}(\cal {L})$.
\end{remark}

\noindent We shall see that the Markov process associated with
$(\cal L',\mu)$-martingale problem is  a subcritical measure-valued
branching process with branching mechanism given by
$$
 \Psi_0(x,z):=\frac{1}{2}\sigma(x)z^2+\int_l^{\infty}\xi\gamma(x,d\xi)z
  +\int_0^l(e^{-z\xi}-1+z\xi)\gamma(x,d\xi).
$$
  For
$i\geq2$, let $\sigma_i:=\sup_x[\int_0^l\xi^i\gamma(x,d\xi)]$. We
first show that each solution of the $(\cal L',\mu)$-martingale
problem has bounded moment of any order.

\begin{lemma}
\label{L2.1}   Suppose that ${\bf Q}'_{\mu}$ is a probability
measure on $D([0,+\infty),M(\mathbb{R}))$ such that under ${\bf
Q}'_{\mu}$ $\omega_0=\mu ~a.s.$ and $\{\omega_t:t\geq0\}$ is a
solution of the $(\cal L',\mu)$-martingale problem. Then for
$n\geq1,t\geq0$ we have
 \begin{eqnarray}
 \label{f2.6}
 {\bf Q}'_{\mu}\{\la1,\omega_t\ra^n\}\ar\leq\ar\sigma_2t/2+
 \la1,\mu\ra^n+C_1(n,\gamma)\int_0^t{\bf Q}'_{\mu}\{\la1,\oz_s\ra\} ds\cr\ar\ar
 +C_2(n,\sigma,\gamma)\int_0^t
 {\bf Q}'_{\mu}\{\la1,\omega_s\ra^{n-1}\}ds
 +C_3(n,\gamma)\int_0^t{\bf Q}'_{\mu}\{\la1,\omega_s\ra^n\}ds,
 \end{eqnarray}
 where $C_1(n,\gamma),~C_2(n,\sigma,\gamma)$ and
 $C_3(n,\gamma)$ are constants which depend on
 $n$, $\sigma$ and $\gamma$.
 \end{lemma}
 \textbf{Proof}.
Let  $n\geq 1$ be fixed. For any $k\geq1$, take $f_k\in
 C_0^2(\mathbb{R})$ such that $f_k(z)=z^n$ for $0\leq
 z\leq k$ and $|f'_k(z)|\leq nz^{n-1}$,
 $f''_k(z)\leq n^2z^{n-2}$ for all $z>k$. Let
 $F_k(\mu)=f_k(\la1,\mu\ra).$ Then $\mathcal{A}F_k(\mu)=0$ and
\begin{eqnarray*}
\mathcal{B}'F_k(\mu)\ar=\ar\frac{1}{2}
 f''_k(\la1,\mu\ra)
 \la\sigma,\mu\ra-\int_{\mathbb{R}}\mu(dx)\int_l^{\infty}\xi
 f'_k(\la1,\mu\ra)\gamma(x,d\xi)\cr
 \ar\ar+\int_{\mathbb{R}}\mu(dx)\int_0^{l}\{
 f_k(\la1,\mu\ra+\xi)-f_n(\la1,\mu\ra)
  -\xi
  f'_k(\la1,\mu\ra)\}\gamma(x,d\xi)\cr
 \ar\leq\ar\frac{1}{2}n^2||\sigma||\la1,\mu\ra^{n-1}
+\sup_x[\int_1^{\infty}\xi\gamma(x,d\xi)]n\la1,\mu\ra^n\cr
  \ar\ar+\int_{\mathbb{R}}\mu(dx)\int_0^{l}
  \frac{1}{2}n^2(\la1,\mu\ra+\xi)^{n-2}\xi^2\gamma(x,d\xi).
 \end{eqnarray*}
 Then we deduce that $$\mathcal{B}'F_k(\mu) \leq
 C_1(n,\gamma)\la1,\mu\ra+C_2(n,\sigma,\gamma)\la1,\mu\ra^{n-1}+
 n\sup_x[\int_1^{\infty}\xi\gamma(x,d\xi)]\la1,\mu\ra^n,$$
where
$C_2(n,\sigma,\gamma)=n^2||\sigma||/2
+\frac{1}{2}\sigma_2n^22^{(n-3)\vee0}$
and
$$C_1(n,\gamma)=
\begin{cases}
n^22^{(n-3)\vee0}\sigma_n/2,&n\geq2,\\
 0,&n=1.
\end{cases}
$$
We have used the Taylor's expansion and elementary inequality
$$(c+d)^{\beta}\leq2^{(\beta-1)\vee0}(c^{\beta}+d^{\beta}),
~\textrm{for
all}~\beta,c,d\geq0.$$ Note that $F_k\in {\cal D}_0(\cal L)$. Thus
$$F_k(\omega_t)-F_k(\omega_0)
-\int_0^t\mathcal{L}'F_k(\omega_s)ds,~t\geq0,$$
is a martingale. We get
\begin{eqnarray*}
 {\bf Q}'_{\mu}f_k(\la1,\omega_t\ra)
 \ar\leq\ar f_k(\la1,\mu\ra)+C_1(n,\gamma)\int_0^t
 {\bf Q}'_{\mu}(\la1,\omega_s\ra)ds\cr
&&+C_2(n,\sigma,\gamma)\int_0^t{\bf
Q}'_{\mu}(\la1,\omega_s\ra^{n-1})ds +C_3(n,\gamma)\int_0^t{\bf
Q}'_{\mu}(\la1,\omega_s\ra^n)ds,
\end{eqnarray*}
where $C_3(n,\gamma)=n\sup_x[\int_1^{\infty}\xi\gamma(x,d\xi)].$ Now
inequality ($\ref{f2.6}$) follows from Fatou's Lemma. \qed
\smallskip

\noindent Observe that, if $F_{m,f}(\mu)=\langle f,\mu^m\rangle$ for
$f\in {C}^{2}(\mathbb{R}^m)$, then
 \begin{eqnarray}
 \label{f2.7}\mathcal
 {A}F_{m,f}(\mu)\ar=\ar\frac{1}{2}\int_{\mathbb{R}^m}\sum_{i=1}^m
 a(x_i)f^{ii}(x_1,\cdots,x_m)\mu^m(dx_1,\cdots,dx_m)\cr
 \ar\ar+\frac{1}{2}\int_{\mathbb{R}^m}\sum_{i,j=1,i\neq j}^m
   \rho(x_i-x_j)f^{ij}(x_1,\cdots,x_m)\mu^m(dx_1,\cdots,dx_m)\cr
 \ar=\ar F_{m,G^mf}(\mu),
 \end{eqnarray}
and
 \begin{eqnarray}
 \label{f2.8}\mathcal {B}'F_{m,f}(\mu)
 \ar=\ar\frac{1}{2}\sum_{i,j=1,i\neq
   j}^m\int_{\mathbb{R}^{m-1}}
   \Psi_{ij}f(x_1,\cdots,x_{m-1})\mu^{m-1}(dx_1,\cdots,dx_{m-1})\cr
 \ar\ar+\sum_{a=2}^m\int_{\mathbb{R}^{m-a+1}}
   \sum_{\{a\}}\Phi_{i_1,\cdots,i_a}f(x_1,\cdots,x_{m-a+1})\mu^{m-a+1}
   (dx_1,\cdots,dx_{m-a+1})\cr
 \ar\ar-\sum_{i=1}^m\int_{\mathbb{R}^m}\int_l^{\infty}\xi
 \gamma(x_i,d\xi)f(x_1,\cdots,x_m)\mu^m(dx_1,\cdots,dx_m),
 \end{eqnarray}
where $\{ a\}=\{1\leq i_1<i_2<\cdots<i_a\leq m\}$. $\Psi_{ij}$
denotes the operator from $B(\mathbb{R}^{m})$ to
$B(\mathbb{R}^{m-1})$ defined by
 \bgeqn
 \label{f2.9}
 \Psi_{ij}f(x_1,\cdots,x_{m-1})=\sigma(x_{m-1})
 f(x_1,\cdots,x_{m-1},\cdots,x_{m-1},\cdots,x_{m-2}),
 \edeqn
where $x_{m-1}$ is in the places of the $i$th and the $j$th
variables of $f$ on the right hand side and $\Phi_{i_1,\cdots,i_a}$
denotes the operator from $B(\mathbb{R}^m)$ to
$B(\mathbb{R}^{m-a+1})$ defined by
 \bgeqn
 \label{f2.10}
 \Phi_{i_i,\cdots,i_a}f(x_1,\cdots,x_{m-a+1})=
 f(x_1,\cdots,x_{m-a+1},\cdots,x_{m-a+1},\cdots,x_{m-a})\int_0^l
 \xi^a\gamma(x_{m-a+1},d\xi),
 \edeqn
where $x_{m-a+1}$ is in the places of the $i_1$th, $i_2$th,
$\cdots$, $i_a$th variables of $f$ on the right hand side. For
$x=(x_1,\cdots,x_m)\in\mathbb{R}^m$, let
$b(x)=\sum_{i=1}^m\int_l^{\infty}\xi\gamma(x_i,d\xi)$. It follows
that
 \begin{eqnarray}
 \label{f2.11}
 \mathcal
 {L}'F_{m,f}(\mu)\ar=\ar F_{m,G^mf}(\mu)-F_{m,bf}(\mu)\cr
 \ar\ar+\frac{1}{2}\sum_{i,j=1,i\neq
 j}^mF_{m-1,\Psi_{ij}f}(\mu)+\sum_{a=2}^m\sum_{\{a\}}
 F_{m-a+1,\Phi_{i_1,\cdots,i_a}f}(\mu).
 \end{eqnarray}
\begin{lemma}
\label{L2.2} Suppose that ${\bf Q}'$ is a probability measure on
$D([0,+\infty),M(\mathbb{R}))$ such that under ${\bf Q}'$
$\{\omega_t:t\geq0\}$ is a solution of the $(\cal
L',\mu)$-martingale problem. Then
 \bgeqn
 \label{f2.12}
 F(\omega_t)-F(\omega_0)-\int_0^t\mathcal{L}'F(\omega_s)ds,~t\geq0,
 \edeqn
under ${\bf Q}'$ is a martingale for each $F(\mu)=F_{m,f}(\mu)=\la
f,\mu^m\ra$ with $f\in C^{2}(\mathbb{R}^m)$.
\end{lemma}
\textbf{Proof}.  For any $k\geq1$, take $f_k\in C_0^2(\mathbb{R}^m)$
such that for $0\leq x_i^2\leq k,~1\leq i\leq m$,
$$f_k(x_1,\cdots,x_m)=\prod_{i=1}^mx_i.$$
For $\{\phi_i\}\subset C^2(\mathbb{R})^{++}$, let
$F_k(\mu)=f_k(\la\phi_1,\mu\ra,\cdots,\la\phi_m,\mu\ra).$ Then
$\lim_{k\rightarrow\infty}F_k(\mu)=F_{m,f}(\mu)$ for all $\mu\in
M({\mathbb{R}})$ and if for every $1\leq i\leq m$,
$0\leq\la\phi_i,\mu\ra^2+l^2||\phi_i||^2\leq k$, we have
$$\mathcal{L}'F_k(\mu)=\mathcal{L}'F_{m,f}(\mu).$$ Introduce a
sequence  stopping times
$$\tau_k:=\inf\{t\geq0,~\textrm{there exists}~i\in\{1,\cdots
m\}~\textrm{such that}~\la\phi_i,\omega_t\ra^2+l^2||\phi_i||^2\geq
k\}\wedge k.$$ Then $\tau_k\rightarrow\infty$ as
$k\rightarrow\infty$. Suppose that $\{H_i\}_{i=1}^n\subset
C(M(\mathbb{R}))$ and $0\leq
 t_1<\cdots<t_n<t_{n+1}$. By Lemma $\ref{L2.1}$ and the
 dominated convergence theorem we deduce that
 \begin{eqnarray*}
 \ar\ar{\bf Q}'\bigg{\{}\big{[}F_{m,f}(\omega_{t_{n+1}})
 -F_{m,f}(\omega_{t_n})
  -\int_{t_n}^{t_{n+1}}\mathcal{L}'F_{m,f}(\omega_s)ds\big{]}
  \prod_{i=1}^nH(\omega_{t_i})\bigg{\}}\cr
\ar=\ar\lim_{k\rightarrow\infty}{\bf Q}'\bigg{\{}
\big{[}F_k(\omega_{t_{n+1}})-F_k(\omega_{t_n})-\int_{t_n}^{t_{n+1}}
 \mathcal{L}'F_k(\omega_s)ds\big{]}\prod_{i=1}^nH(\omega_{t_i})\bigg{\}}\cr
\ar\ar+\lim_{k\rightarrow\infty}{\bf
Q}'\bigg{\{}\big{[}\int_{t_n}^{t_{n+1}}
 \mathcal{L}'F_k(\omega_s)1_{\{\tau_k\leq s\}}ds-\int_{t_n}^{t_{n+1}}
 \mathcal{L}'F_{m,f}(\omega_s)1_{\{\tau_k\leq
 s\}}ds\big{]}\prod_{i=1}^nH(\omega_{t_i})\bigg{\}}\cr
\ar=\ar0.
 \end{eqnarray*}
That is under ${\bf Q}'$
$$F_{m,f}(\omega_{t})-F_{m,f}(\omega_{0})
  -\int_{0}^{t}\mathcal{L}'F_{m,f}(\omega_s)ds,~t\geq0,
  $$
 is a martingale  for $f=\prod_{i=1}^m\phi_i$
with $\{\phi_i\}\subset C^2(\mathbb{R})^{++}$(and therefore
$\{\phi_i\}\subset C^2(\mathbb{R})$). Since $f\in C^2(\mathbb{R})$
can be approximated by polynomials in such a way that not only $f$
but its derivatives up to second order are approximated uniformly on
compact sets, by an approximating procedure (\ref{f2.12}) is a
martingale for $F(\mu)=\la f,\mu^m\ra$ with $f\in C_c^2(\mbb R^m)$
(see \cite{[EK86]}, p.501). By Remark $\ref{rem1.1}$, (1,~0) is in
the bp-closure of $G^m$. In fact, let $h\in C_c^2(\mbb R^m)$ satisfy
$1_{B(0,1)}\leq h\leq 1_{B(0,2)}$ and $h_k(x)=h(x/k)\in C_c^2(\mbb
R^m)$. Then for $f\in C^2(\mbb R^m)$, we can approximate $(f,G^mf)$
by $\{fh_k,G^mfh_k\}$. According to  ($\ref{f2.11}$) and Lemma
\ref{L2.1}, we see the desired result follows by another
approximating procedure. \qed

\noindent Let $G_b^m:=G^m-b$. By Theorem 5.11 of \cite{[Dy65]},
there exists a diffusion process on $\hat{C}(\mbb{R}^m)$ generated
by $G_b^m\mid_{C_c^2(\mbb{R}^m)}$ (and therefore
$G_b^m\mid_{\hat{C}^2(\mbb{R}^m)}$). Its transition density
$q^m(t,x,y)$ is the fundamental solution of the equation
 \bgeqn
 \label{f2.13}
 \frac{\partial u}{\partial t}=G_b^mu.
 \edeqn
  The semigroup corresponding to the operator $G_b^m$ is defined by
  \bgeqn
  \label{f2.14}
   T_t^mf(x)=\int q^m(t,x,y)f(y)dy
  \edeqn
   for $f\in \hat{C}(\mbb{R}^m)$ and  can therefore be extended
   to all of $B(\mathbb{R}^m)$.
  According to $0.24.A_2$ of \cite{[Dy65]}, for $f\in C(\mbb R^m)$
  $$
  \lim_{t\rightarrow0}\int q^m(t,x,y)f(y)dy=f(x)~~~(x\in \mbb R^m),
  $$
  where the convergence is uniform on every bounded subset.
   On the other hand,  $(T_t^m)_{t\geq0}$ is strong Feller, i.e.,
    for $f\in  B(\mbb R^m)$ and $t>0$,
   $T_t^mf\in C(\mbb R^m)$. In fact, according to $1^{\circ}$ of the
   proof of Theorem 5.11 of \cite{[Dy65]}, $T^m_t f\in C^2(\mbb
   R^m)$ satisfies equation (\ref{f2.13}). Hence for $f\in C^2(\mbb
   R^m)$
   $$
   \frac{T_t^mf(x)-f(x)}{t}=\frac{u_t(x)-f(x)}{t}=\frac{1}{t}
   \int_0^tG_b^mu_s(x)ds.
   $$
 Therefore
 $$
 \lim_{t\rightarrow0}\frac{T_t^mf(x)-f(x)}{t}=G_b^mf(x),
 $$
 where the convergence is bounded and pointwise.
Let $\tilde{G}_b^m$ denote the weak generator of $(T_t^m)_{t\geq0}$.
Thus $C^2(\mbb R^m)$ belong to the domain of $\tilde{G}_b^m$ and
$T_t^m C^2(\mbb R^m)\subset C^2(\mbb R^m)$. Also,
$\tilde{G}_b^m|_{C^2(\mbb R^m)}=G^m_b|_{C^2(\mbb R^m)}$. Let
$p^m(t,x,y)$ denote the transition density corresponding to the
semigroup $(P_t^m)_{t\geq0}$. According to $6^\circ$ of the proof of
Theorem 5.11 of \cite{[Dy65]}, we see for all $t>0$, $x\in\mbb R^m$,
$A\in \mathfrak{B}(\mbb R^m)$,
$$
\int_{A}p^m(t,x,y)dy\geq\int_{A}q^m(t,x,y)dy.
$$
Therefore, for $f\in B(\mbb R^m)^+$, $$P_t^mf(x)\geq T_t^mf(x).$$
Next, we  define a dual process and  reveal its connection to the
solutions of the $({\cal L}',\mu)$-martingale problem.

\medskip
\noindent Let $\{M_t:t\geq0\}$ be a nonnegative integer-valued
c\`{a}dl\`{a}g Markov process. For  $i\geq j$, the transition
intensities $\{q_{ij}\}$ defined by
$$
 q_{ij}=\left\{\begin{array}{ccc}
                 \sum_{i\neq j}-q_{ij} &if&j=i  \cr
                 \frac{1}{2}i(i-1)+\left
                 (\begin{array}{c}i\cr2\end{array}\right)&if&j=i-1\cr
                 \left(\begin{array}{c}i\cr j-1\end{array}\right)
                               &if&1\leq j\leq i-2\cr
                 \end{array}
                 \right.
$$
and $q_{ij}=0$ for $i<j$. Let $\tau_0=0$ and $\tau_{M_0}=\infty$,
and let $\{\tau_k:1\leq k\leq M_0-1\}$ be the sequence of jump times
of $\{M_t:t\geq0\}$. That is $\tau_1=\inf\{t\geq0:M_t\neq
M_0\},\cdots,\tau_k=\inf\{t>\tau_{k-1}:M_t\neq M_{\tau_{k-1}}\} .$\\
Let $\{\Gamma_k:1\leq k\leq M_0-1\}$ be a sequence of random
operators which are conditionally independent given $\{M_t :t\geq
0\}$ and satisfy
$$\textbf{P}\{\Gamma_k=\Psi_{ij}|M(\tau_k-)=l,M(\tau_k)=l-1\}
=\frac{1}{2l(l-1)} ,\textrm{\ \ \ }1\leq i\neq j\leq l,$$
$$\textbf{P}\{\Gamma_k=\Phi_{i_1,i_2}|M(\tau_k-)=l,M(\tau_k)=l-1\}=
\frac{1}{l(l-1)},\textrm{\ \ \ }1\leq i_1<i_2\leq l,$$ and for
$a\geq3$,
$$\textbf{P}\{\Gamma_k=\Phi_{i_1,\cdots,i_a}|M(\tau_k-)=l,M(\tau_k)
=l-a+1\}=\frac{1}{\left(\begin{array}{c}l\cr a\end{array}\right)
},\textrm{\ \ \ }1\leq i_1<\cdots<i_a\leq l,$$ where $\Psi_{ij}$ and
$\Phi_{i_1,\cdots,i_a}$ are defined by (\ref{f2.9}) and
(\ref{f2.10}) respectively. Let $\textbf{B}$ denote the topological
union of $\{B(\mathbb{R}^m):m = 1,2,\cdots\}$ endowed with pointwise
convergence on each $B(\mathbb{R}^m)$. Then
 \bgeqn
 \label{f2.15}Y_t={T}_{t-\tau_k}^{M_{\tau_k}}\Gamma_k
 {T}_{\tau_k-\tau_{k-1}}^{M_{\tau_{k-1}}}\Gamma_{k-1}\cdots
 {T}_{\tau_2-\tau_1}^{M_{\tau_1}}\Gamma_1{T}_{\tau_1}^{M_0}Y_0, \textrm{\
 \ \ }\tau_k\leq t<\tau_{k+1},~~0\leq k\leq M_0-1,
 \edeqn
defines a Markov process $\{Y_t:t\geq0\}$ taking values from ${\bf
B}$. Clearly, $\{(M_t,Y_t):t\geq0\}$ is also a Markov process. Let
$\textbf{E}_{m,f}^{\sigma,\gamma}$ denote the expectation given
$M_0=m$ and $Y_0=f\in B(\mathbb{R}^m)$.
\begin{theorem}
\label{T2.1}
  Suppose that $\{X_t:t\geq0\}$ is
a c\`{a}dl\`{a}g $M({\mathbb{R}})$-valued process. If $\{X_t:t\geq
0\}$ is a solution of the $({\mathcal{L}}',\mu)$-martingale problem
and assume that $\{X_t:t\geq0\}$ and $\{(M_t,Y_t):t\geq0\}$ are
defined on the same probability space and independent of each other,
then
 \bgeqn
 \label{f2.16}
 \textbf{E}\left\langle f,X_t^m\right\rangle
 =\mathbf{E}_{m,f}^{\sigma,\gamma}\big{[}\left\langle
 Y_t, \mu^{M_t}\right\rangle
 \exp\big{\{}\int_0^t(2^{M_s}
 +\frac{M_s(M_s-1)}{2}-M_s-1)ds\big{\}}\big{]}
 \edeqn
for any $t\geq0$, $f\in B(\mathbb{R}^m)$ and integer $m\geq1$.
\end{theorem}
\textbf{Proof}. In this proof we set
$F_{\mu}(m,f)=F_{m,f}(\mu)=\langle f,\mu^m \rangle$. By Lemma
$\ref{L2.1}$, we have that for each $m\geq1$,
$\mathbf{E}[\la1,X_t\ra^m]$ is a locally bounded function of
$t\geq0$. Then by  martingale inequality we have that
$\mathbf{E}[\sup_{0\leq s\leq t}\la1,X_s\ra^m]$ is a locally bounded
function of $t\geq0$.

\noindent By the definition of $Y$ and elementary properties of $M$,
we know that $\{(M_t,Y_t): t\geq0\}$ has weak generator $\mathcal
{L}^*$ given by
 \begin{eqnarray}
 \label{f2.17}
 \mathcal {L}^*F_{\mu}(m,f)\ar=
 \ar F_{\mu}(m,{G}_b^mf)+
   \frac{1}{2}\sum_{i,j=1,i\neq
   j}^m\left[F_{\mu}(m-1,\Psi_{ij}f)-F_{\mu}(m,f)\right]\cr
 \ar\ar+\sum_{k=2}^m\left(\sum_{\{1\leq i_1<\cdots<i_k\leq m\}}
 \left[F_{\mu}(m-k+1,\Phi_{i_1,\cdots,\i_k}f)
 -F_{\mu}(m,f)\right]\right)\end{eqnarray}
with $f\in C^2(\mbb{R}^m)$. In view of (\ref{f2.11}) we have
 \bgeqn
 \label{f2.18}
 \mathcal {L}^*F_{\mu}(m,f)={\mathcal
 {L}}'F_{m,f}(\mu)-(2^m+\frac{1}{2}m(m-1)-m-1)F_{\mu}(m,f).
 \edeqn
Then it is easy to verify that the inequalities in Theorem 4.4.11 of
\cite{[EK86]} are satisfied. Then the desired conclusion follows
from Corollary 4.4.13 of \cite{[EK86]}. \qed

\begin{remark}
\label{gap} We point out that there exists a gap in the proof of
establishing the dual relationship of \cite{[DLW01]}. There it was
assumed $\sigma$ is a bounded measurable function and $\gamma=0$.
When they established the dual relationship, they  used a
relationship which is similar to ($\ref{f2.18}$). However, note that
($\ref{f2.18}$)  makes sense if $f\in\mathcal{D}({G}_b^m)$ and $Y_t$
need not always take values in $\mathcal{D}({G}_b^m)$ if we only
assume that $\sigma$ is a bounded measurable function and $G^m$ is
elliptic. If we assume that $\sigma\in C_{\partial}^2(\mbb{R})$ and
$G^m$ is uniformly elliptic, then the argument there can be applied
to establish the dual relationship there. If $c=0$, $G^m$ need not
always be uniformly elliptic. Our methods cannot be applied to
obtain the uniqueness of the corresponding martingale problem.
Dawson and Li \cite{[DL03]} constructed SDSM from one-dimensional
excursion when $c=0$ and $\gamma(x,d\xi)=0$. From the construction
there,  an important property of the SDSM was revealed. That is when
$c=0$, the process always lives in the space of purely atomic
measures. We can also follow the idea there to construct
discontinuous SDSM.
\end{remark}

\begin{theorem}
\label{T2.2}  Suppose that for each $\mu\in M(\mathbb{R})$ there is
a probability measure ${\bf Q}'_{\mu}$ on
$D([0,\infty),M(\mathbb{R}))$ such that ${\bf
Q}'_{\mu}\{\la1,\omega_t\ra^m\}$ is locally bounded in $t\geq0$ for
every $m\geq1$ and such that $\{\omega_t:t\geq0\}$ under ${\bf
Q}'_{\mu}$ is a solution of the $(\cal L',\mu)$-martingale problem.
Then ${\bf Q}':=\{{\bf Q}'_{\mu}:\mu\in M(\mathbb{R})\}$ defines a
Markov process with transition semigroup $(Q'_t)_{t\geq0}$ given by
 \bgeqn
 \label{f2.19}
 \int_{M(\mathbb{R})}\left\la f,\nu^m\right\ra Q'_t(\mu,d\nu)
 =\mathbf{E}_{m,f}^
 {\sigma,\gamma}\left[\left\langle
 Y_t,\mu^{M_t}\right\rangle
 \exp\left\{\int_0^t\left(2^{M_s}+\frac{M_s(M_s-1)}{2}
 -M_s-1\right)ds\right\}\right]
 \edeqn
 for $f\in B(\mbb R^m)$.
\end{theorem}
\textbf{Proof}. Let $Q'_t(\mu,\cdot)$ denote the distribution of
$\omega_t$ under ${\bf Q}'_{\mu}$. By Theorem $\ref{T2.1}$, we
obtain ($\ref{f2.19}$).  We first consider the case that
$\sigma(x)\equiv\sigma_0$ for a constant $\sigma_0$ and
$\gamma(x,d\xi)\equiv\hat{\gamma}(d\xi)$ such that
$\int_l^{\infty}\hat{\gamma}(d\xi)=0$. In this case,
$\{\la1,\omega_t\ra:t\geq0\}$ is a critical continuous state
branching process with generator $\mathfrak{L}$ given by
 \begin{eqnarray}
\mathfrak{L}f(x)&=&\frac{1}{2}\sigma_0xf''(x)+x\int_0^l
\left(f(x+\xi)-f(x)-\xi f'(x)\right)\hat{\gamma}(d\xi)
\end{eqnarray}
for $f\in C^2(\mathbb{R})$. By Kawazu and Watanabe \cite{[KW71]} we
deduce that
$$\int_{M(\mathbb{R})}e^{\lambda\la1,\nu\ra}Q_t(\mu,d\nu)
=e^{\la1,\mu\ra\varphi(t,\lambda)},~~t\geq0,~\lambda\geq0,$$ where
$\varphi(t,\lambda)$ is the solution of
$$\begin{cases}\frac{\partial\varphi}{\partial
t}(t,\lz)=R(\varphi(t,\lz)),\\
\varphi(0,\lz)=\lz, \end{cases}$$  and $R(\lz)$ is given as follows:
$$R(\lz)=-\frac{1}{2}\sigma_0\lz^2
-\int_0^{l}(e^{-\lz\xi}-1+\lz\xi)\hat{\gamma}(d\xi).$$ Then for each
$f\in B(\mathbb{R})^+$ the power series
 \bgeqn
 \label{f2.21}
 \sum_{m=0}^{\infty}\frac{1}{m!}\int_{M(\mathbb{R})}\la
 f,\nu\ra^mQ'_t(\mu,d\nu)\lz^m
 \edeqn
has a positive radius of convergence. By this and Theorem 30.1 of
\cite{[B86]}, it is easy to show that $Q'_t(\nu,\cdot)$ is the
unique probability measure on $M(\mathbb{R})$ satisfying
($\ref{f2.19}$). Now the result follows from Theorem 4.4.2 of
\cite{[EK86]}. For general case, let $\sigma_0=||\sigma||$ and
$f^{\otimes m}(x_1,\cdots,x_m)=f(x_1)\cdots f(x_m)$. We can find a
measure $\hat{\gamma}(d\xi)$ on $(0,+\infty)$ such that for every
$k\geq2$
$$
C_{\gamma}:=\sup_x[\int_0^1\xi^2\gamma(x,d\xi)
+\int_1^l\xi\gamma(x,d\xi)]
\leq\int_0^l\xi^k\hat{\gamma}(d\xi)<\infty
$$
and $\int_l^{\infty}\hat{\gamma}(d\xi)=0$. In fact, since $l>1$, we
can let $\hat{\gamma}(d\xi)=(k_l+1)C_{\gamma}1_{(0,l)}(\xi)d\xi$,
where $d\xi$ denotes the Lebesgue measure and
$k_l=\min\{k\geq2:l^k/(k+1)>1\}$. We obtain that for each $k\geq2$
$$
\sup_x[\int_0^l\xi^k\gamma(x,d\xi)]\leq
l^k\int_0^l\xi^k\hat{\gamma}(d\xi).
$$
By ($\ref{f2.19}$) and ($\ref{f2.15}$) we have
$$
\int_{M(\mathbb{R})}\la f,\nu\ra^m Q'_t(\mu,d\nu)\leq\mathbf{E}_{m,
l^mf^{\otimes m}}^
 {\sigma_0, \hat{\gamma}}\left[\left\langle
 Y_t,\mu^{M_t}\right\rangle
 \exp\left\{\int_0^t\left(2^{M_s}+\frac{M_s(M_s-1)}{2}
 -M_s-1\right)ds\right\}\right]
 $$
 for $f\in B(\mbb R)^+.$
Then the power series ($\ref{f2.21}$) also has a positive radius of
convergence and the desired result follows as in previous case. \qed
\begin{remark}
From (\ref{f2.11}), we may regard the Markov process associated with
$({\cal L}',\mu)$-martingale problem as a measure-valued branching
process with branching mechanism given by
$$\Psi_1({x,z}):=\frac{1}{2}\sigma(x)z^2
+\int_0^l(e^{-z\xi}-1+z\xi)\gamma(x,d\xi).$$
 and its spatial motion is a diffusion process generated by
$$
\frac{1}{2}\sum_{i=1}^ma(x_i)\frac{\partial^2}{\partial x_i^2}+
\frac{1}{2}\sum_{i,j=1,i\neq
j}^m\rho(x_i-x_j)\frac{\partial^2}{\partial x_i\partial
x_j}-\sum_{i=1}^m\int_l^{\infty}\xi\gamma(x_i,d\xi)
$$
which represents  `$G^m$-diffusion  killed at a rate
$\sum_{i=1}^m\int_l^{\infty}\xi\gamma(x_i,d\xi)$'; see Rogers and
Williams \cite{[RW00]} and references therein for more details of
`Markov process with killing'.
\end{remark}

\subsection{Uniqueness for $(\mathcal{L},\mu)$-martingale problem}

In this section, we will consider a localization procedure suggested
by Stroock \cite{[S75]} to show that the uniqueness for the
$(\mathcal{L}, \mu)$-martingale problem follows from the uniqueness
of the $(\mathcal{L}',\mu)$-martingale problem. Although some
arguments in this subsection are similar to those of  \cite{[ER91]}
and \cite{[S75]}, we shall give the details for the convenience of
the reader. We assume that the for each $\mu\in M(\mbb R)$, $(\cal
L',\mu)$-martingale problem is well-posed. The existence for the
$(\cal L', \mu)$-martingale problem will be revealed in Section 3.
Let  $\bf Q'$ denote the Markovian system defined in Theorem
$\ref{T2.2}$. Let ${\bf Q}'_{s,\mu}={\bf Q}'(\cdot|\omega_s=\mu)$.
Then ${\bf Q}'_{s,\mu}$ is  also a Markovian system starting from
$(s,\mu)$ whose transition semigroup is the same with $\bf Q'$.

\noindent Let $\{\omega_t: t\geq0 \}$ denote the coordinate process
of $D([0,\infty),M(\mathbb{R}))$. Let
$\Omega=D([0,\infty),M(\mathbb{R}))$. Set
$\mathcal{F}_t=\sigma\{\omega_s:0\leq s\leq t\}$, and take
$\mathcal{F}^t=\sigma\{\oz_s: t\leq s\}$.
\begin{definition}
For $\mu\in M(\mathbb{R})$, we say a probability measure
$\mathbf{Q}_{s,\mu}$ on $(\Omega, \mathcal{F}^s)$ is a solution of
the $(\mathcal {L},\mu)$-\textsl{martingale problem} if
$\mathbf{Q}_{s, \mu}(\omega_s=\mu)=1$ and
 \bgeqn \label{f2.22}
  F(\omega_t)-F(\mu)-\int_s^t\mathcal {L}F(\omega_u)du,
  \quad t\geq s,
 \edeqn
 is a local martingale for each $F\in\mathcal {D}(\mathcal {L})$.
\end{definition}
In the following we will write $\mathbf{Q}_{\mu}$ instead of
$\mathbf{Q}_{0,\mu}$ and write $\cal F$ instead of ${\cal F}^{0}$.
Let $S({\mathbb{R}})$ denote the space of finite signed Borel
 measures on ${\mathbb{R}}$ endowed with the $\sigma$-algebra
 generated by the mappings $\mu\mapsto\la f,\mu\ra$ for all $f\in
 C({\mathbb{R}})$. Let
 $S({\mathbb{R}})^{\circ}=S({\mathbb{R}})\setminus\{0\}$
 and $M({\mathbb{R}})^{\circ}=M({\mathbb{R}})\setminus\{0\}$.
 The
 following theorem is analogous to  Th\'{e}or\`{e}m 7 of \cite{[ER91]}.
\begin{theorem}
\label{T2.3} Suppose that  a probability measure $\mathbf{Q}_{\mu}$
on $(\Omega, \mathcal{F})$ is a solution of the $(\cal
L,\mu)$-martingale problem. Define an optional random measure
 $N(ds,d\nu)$ on $[0,\infty)\times S({\mathbb{R}})^{\circ}$ by
 $$
 N(ds,d\nu)=\sum_{s>0}1_{\{\Delta
 \omega_s\neq0\}}\delta_{(s,\Delta \omega_s)}(ds,d\nu),
 $$
 where $\Delta \oz_s=\oz_s-\oz_{s-}\in S({\mathbb{R}})$. Let
 $\hat{N}(ds,d\nu)$ denote the predictable compensator of $N(ds,d\nu)$
 and let $\tilde{N}(ds,d\nu)$ denote the corresponding martingale measure
  under ${\bf Q}_{\mu}$.
 Then $\hat{N}(ds,d\nu)=dsK(\oz_s,d\nu)$ with $K(\mu,d\nu)$
given by
$$
\int_{M({\mathbb{R}})^{\circ}}F(\nu)K(\mu,d\nu)=
 \int_{\mbb R}\mu(dx)\int_0^{\infty}F(\xi\dz_x)\gamma(x,d\xi),
$$
and  for $\phi\in C^2(\mbb R)^+$,
 \bgeqn
 \label{f2.23}
 M_t(\phi):=\la\phi,\oz_t\ra-\la\phi,\mu\ra
 -\frac{1}{2}\int_0^t\la a\phi'',\oz_s\ra ds
 ,\quad t\geq0,
 \edeqn
is a martingale and we also have that
$$
M_t(\phi)=M_t^c(\phi)+M_t^d(\phi),
$$
 where  $M_t^c(\phi)$  under $\bf Q_{\mu}$ is a continuous
 martingale with quadratic variation process given by
 \bgeqn
 \label{f2.24}
 \la M^c(\phi)\ra_t=\int_0^t\la\sigma\phi^2,\omega_s\ra ds+
 \int_0^tds\int_{\mathbb{R}}\la h(z-\cdot)\phi',\omega_s\ra^2
 dz,
 \edeqn
 and
 \bgeqn
 \label{f2.25}
 M_t^d(\phi)=\int_0^{t+}
 \int_{M({\mathbb{R}})^{\circ}}\la\phi,\nu\ra\tilde{N}(ds,d\nu)
 \edeqn
 is a purely discontinuous  martingale under $\bf Q_{\mu}$.
\end{theorem}

\smallskip
\noindent$\bf Proof$. Some arguments in the proof of this theorem
are similar to those of Theorem 6.1.3 of \cite{[D93]}. The proof
will be divided into 4 steps.

\textit{Step 1}.
 Since $e^{-\la\phi,\nu\ra}\in {\cal
D}_0(\cal L)$ for $\phi\in C^2(\mbb R)^{++}$,
 \bgeqn
 \label{f2.26}
 W_t(\phi):=e^{-\la\phi,\oz_t\ra}
 -\int_0^t e^{-\la\phi,\oz_s\ra}[-\frac{1}{2}\la
 a\phi'',\oz_s\ra+\frac{1}{2}\int_{{\mathbb{R}}}\la h(z-\cdot)
 \phi',\oz_s\ra^2dz+
 \la \Psi(\phi), \oz_s\ra]ds, ~t\geq0,
 \edeqn
 is a ${\bf Q}_{\mu}$-martingale with $\phi\in C^2(\mathbb{R})^{++}$,
  where $\Psi(\phi):=\Psi(x,\phi(x))$. Therefore, $\{W_t(\phi)\}$ is a
  local martingale for $\phi\in C^2(\mbb R)^+$.
  Let
 $$Z_t(\phi):=\exp\{-\la\phi,\oz_t\ra\},$$
 $$H_t(\phi):=\exp\bigg{\{}-\la\phi,\oz_t\ra+\int_0^t\big{[}\frac{1}{2}\la
 a\phi'',\oz_s\ra-\frac{1}{2}\int_{{\mathbb{R}}}\la
 h(z-\cdot)\phi',\oz_s\ra^2dz-\la\Psi(\phi),\oz_s\ra\big{]}ds\bigg{\}}
 $$
 and
 $$Y_t(\phi):=\exp\bigg{\{}\int_0^t\big{[}\frac{1}{2}\la
 a\phi'',\oz_s\ra-\frac{1}{2}\int_{{\mathbb{R}}}\la
 h(z-\cdot)\phi',\oz_s\ra^2dz-\la\Psi(\phi),\oz_s\ra\big{]} ds\bigg{\}}.
 $$
  By integration by parts,
 \begin{eqnarray*}
 & &\int_0^{t}Y_s(\phi)dW_s(\phi)\cr
 &=&\int_0^{t}Y_s(\phi)dZ_s(\phi)\cr
 & &~-\int_0^{t}Y_s(\phi)e^{-\la\phi,\oz_s\ra}\big{[}-\frac{1}{2}\la
 a\phi'',\oz_s\ra+\frac{1}{2}\int_{{\mathbb{R}}}\la
  h(z-\cdot)\phi',\oz_s\ra^2dz+
 \la \Psi(\phi), \oz_s\ra\big{]}ds
 \cr&=&H_{t}(\phi)-Z_0(\phi)
 \end{eqnarray*}
is a ${\bf Q}_{\mu}$-local martingale. We also have
 $$Z_t(\phi)=Y_t^{-1}(\phi)H_t(\phi),$$
 and, again by integration by parts,
 \begin{eqnarray}
 \label{f2.27}
 dZ_t(\phi)
 \ar=\ar Y_t^{-1}(\phi)dH_t(\phi)+H_{t-}(\phi)dY_t^{-1}(\phi)\cr
 \ar=\ar Y_t^{-1}(\phi)dH_t(\phi)\cr
 \ar\ar+Z_{t-}(\phi)\big{[}-\frac{1}{2}\la
 a\phi'',\oz_{t-}\ra+\frac{1}{2}
 \int_{{\mathbb{R}}}\la h(z-\cdot)\phi',\oz_{t-}\ra^2dz+
 \la \Psi(\phi), \oz_{t-}\ra\big{]}dt.
 \end{eqnarray}
Then $\{Z_t(\phi):t\geq0\}$ is a special semi-martingale with
$\phi\in C^2(\mathbb{R})^{+}$ (see Definitions 1.4.21 of
\cite{[JS87]}).

\textit{Step 2}. By the same argument as in the proof of Lemma
\ref{L2.1}, we have that
$$
{\bf Q}_{\mu}[\oz_t(1)]\leq
\la1,\mu\ra+C_1(\sigma,\gamma)\int_0^t{\bf Q}_{\mu}[\oz_s(1)]ds,
$$
where
$C_1(\sigma,\gamma):=||\sigma||+2\sup_x\int_1^{\infty}\xi\gamma(x,d\xi)+
\sup_x\int_0^1\xi^2\gamma(x,d\xi)$. By Gronwall's inequality

 \bgeqn \label{f2.28}{\bf Q}_{\mu}[\oz_t(1)]\leq \la1,\mu\ra
e^{C_1(\sigma,\gamma)t}.
 \edeqn
For any $k\geq1$, take $f_k\in C_0^2(\mbb R)$ such that $f_k(x)=x$
for $|x|\leq k$ and $|f_k^\prime(x)|\le 1$ for all $x\in \mbb{R}$.
We see for each $\phi\in C^2(\mbb R)^{++}$,
$$
\lim_{k\rightarrow\infty}f_k(\la\phi,\mu\ra)=\la\phi,\mu\ra\quad\textrm{and}
\quad \lim_{k\rightarrow\infty}{\cal
L}f_k(\la\phi,\mu\ra)=\frac{1}{2}\la a\phi'',\mu\ra.
$$
Since $f_k(\la\phi,\mu\ra)\in {\cal D}_0(\cal L)$, by (\ref{f2.28})
and dominated convergence theorem an approximation argument shows
that for $\phi\in C^2(\mbb R)^{++}$
$$
\la\phi,\oz_{t}\ra=\la\phi,\mu\ra+\frac{1}{2}\int_0^{t}\la
a\phi'',\oz_s\ra ds
 +M_{t}(\phi),
$$
 where $\{M_{t}(\phi):t\geq0\}$ is a  martingale.
For $\phi\in C^2(\mbb R)^{+}$, we have $\{M_{t}(\phi+\varepsilon)\}$
are  martingales for $\varepsilon>0$. By letting
$\varepsilon\rightarrow0$, (\ref{f2.28}) ensures that
 $$
 M_{t}(\phi)=\la\phi,\oz_{t}\ra-\la\phi,\mu\ra-\frac{1}{2}\int_0^{t}\la
a\phi'',\oz_s\ra ds,\quad t\geq0,
 $$
  is a martingale for $\phi\in C^2(\mbb R)^{+}$.
 By Corollary 2.2.38 of  \cite{[JS87]}, $\{M_t(\phi)\}$ admits
 a unique representation
 $$
 M_t(\phi)=M_t^c(\phi)+M_t^d(\phi),
 $$
 where $\{M_t^c(\phi)\}$ is a continuous local martingale with
quadratic variation process $\{C_t(\phi)\}$ and
 \bgeqn
 \label{f2.29}
 M_t^d(\phi)=\int_0^{t+}\int_{S({\mathbb{R}})^{\circ}}
 \la\phi,\nu\ra\tilde{N}(ds,d\nu)
 \edeqn
is a purely discontinuous local martingale.  Moreover,
$\{\la\phi,\oz_t\ra\}$ is a semimartingale. An application of
It\^{o}'s formula for semimartingale (see Theorem 1.4.57 of
\cite{[JS87]}) yields
 \beqlb
 \label{f2.30}
 dZ_t(\phi)&=&Z_{t-}(\phi)[-dU_t(\phi)+\frac{1}{2}dC_t(\phi)
 +\int_{S({\mathbb{R}})^{\circ}}
 (e^{-\la\phi,\nu\ra}-1
 +\la\phi,\nu\ra)N(dt,d\nu)]\cr
 & &+d(loc.mart.),
 \eeqlb
where $U_t(\phi)=\frac{1}{2}\int_0^t\la a\phi'',\oz_s\ra ds$ is of
locally bounded variation. Note that
$$
0\leq Z_{s-}(\phi)(e^{-\la\phi,\nu\ra}-1
 +\la\phi,\nu\ra)\leq C(|\la \phi,\nu\ra|\wz|\la\phi,\nu\ra^2|)
$$
for some constant $C\geq0$. According to Theorem 1.4.47 of
\cite{[JS87]}, $\sum_{s\leq t}(\la\phi,\Delta \oz_s\ra)^2<\infty$.
Thus the first term in (\ref{f2.30}) has finite variation over each
finite interval $[0,t]$.  Since $\{Z_t(\phi)\}$ is a special
semimartingale, Proposition 1.4.23 of \cite{[JS87]} implies that
$$\int_0^{t+}\int_{S({\mathbb{R}})^{\circ}}
 Z_{s-}(\phi)(e^{-\la\phi,\nu\ra}-1
 +\la\phi,\nu\ra){N}(ds,d\nu)$$
 is of locally integrable variation. Thus it is locally integrable.
According to Proposition 2.1.28 of \cite{[JS87]},
\begin{eqnarray*}
\ar\ar\int_0^{t+}\int_{S({\mathbb{R}})^{\circ}}
 Z_{s-}(\phi)(e^{-\la\phi,\nu\ra}-1
 +\la\phi,\nu\ra)\tilde{N}(ds,d\nu)\cr\ar
 =\ar\int_0^{t+}\int_{S({\mathbb{R}})^{\circ}}
 Z_{s-}(\phi)(e^{-\la\phi,\nu\ra}-1
 +\la\phi,\nu\ra){N}(ds,d\nu)\cr\ar\ar-\int_0^{t+}\int_{S({\mathbb{R}})^{\circ}}
 Z_{s-}(\phi)(e^{-\la\phi,\nu\ra}-1
 +\la\phi,\nu\ra)\hat{N}(ds,d\nu)
\end{eqnarray*}
is a purely discontinuous local martingale. Therefore,
 \beqlb
 \label{f2.31}
 dZ_t(\phi)&=&Z_{t-}(\phi)[-dU_t(\phi)+\frac{1}{2}dC_t(\phi)
 +\int_{S({\mathbb{R}})^{\circ}}
 (e^{-\la\phi,\nu\ra}-1
 +\la\phi,\nu\ra)\hat{N}(dt,d\nu)]\cr
 & &+d(loc.mart.).
 \eeqlb

\textit{Step 3}.
   Since $Z_t(\phi)$ is a special
semimartingale we can identify the predictable components of locally
integrable variation in the two decompositions ($\ref{f2.27}$) and
($\ref{f2.31}$) to get that
 \begin{eqnarray*}
 &&Z_{t-}(\phi)[-\frac{1}{2}\la
 a\phi'',\oz_{t-}\ra+\frac{1}{2}\int_{{\mathbb{R}}}\la
 h(z-\cdot)\phi',\oz_{t-}\ra^2dz+
 \la \Psi(\phi), \oz_{t-}\ra]dt\cr
 &=&Z_{t-}(\phi)[-dU_t(\phi)+\frac{1}{2}dC_t(\phi)
 +\int_{S({\mathbb{R}})^{\circ}}
         (e^{-\la\phi,\nu\ra}-1
         +\la\phi,\nu\ra)\hat{N}(dt,d\nu)].
 \end{eqnarray*}
Then
 \beqlb
 \label{f2.32}
 &&\int_0^t[-\frac{1}{2}\la
 a\phi'',\omega_{s}\ra+\frac{1}{2}
 \int_{{\mathbb{R}}}\la h(z-\cdot)\phi',\oz_{s}\ra^2dz+
 \la \Psi(\phi), \oz_{s}\ra]ds\cr
 &=&-U_t(\phi)+\frac{1}{2}C_t(\phi)
 +\int_0^t\int_{S({\mathbb{R}})^{\circ}}
 (e^{-\la\phi,\nu\ra}-1+\la\phi,\nu\ra)\hat{N}(ds,d\nu).
 \eeqlb
 According to  ($\ref{f2.28}$) and ($\ref{f2.29}$), we can deduce
 that  $C_t(\theta\phi)=\theta^2C_t(\phi)$ with $\theta>0$.
  Replacing
 $\phi$ by $\theta\phi$ with $\theta>0$ in ($\ref{f2.32}$), we have
\beqlb
 \label{f2.33}
 &&-\theta\int_0^t\frac{1}{2}\la
 a\phi'',\oz_{s}\ra ds+\frac{\theta^2}{2}\int_0^t
 \int_{{\mathbb{R}}}\la h(z-\cdot)\phi',\oz_{s}\ra^2dzds+
 \frac{\theta^2}{2}\int_0^t\la\sigma\phi^2,\oz_s\ra ds\cr
 &&
 +\int_0^tds\int_{{\mathbb{R}}}\oz_s(dx)
 \int_0^{\infty}\gamma(x,d\xi)(e^{-\theta\xi\phi(x)}-1+\theta\xi\phi(x))\cr
 &=&-\theta U_t(\phi)+\frac{\theta^2}{2}C_t(\phi)
 +\int_0^t\int_{S({\mathbb{R}})^{\circ}}
 (e^{-\theta\la\phi,\nu\ra}-1+\theta\la\phi,\nu\ra)\hat{N}(ds,d\nu).
 \eeqlb
We conclude that
 \bgeqn
\label{f2.34}C_t(\phi)=\int_0^tds\int_{{\mathbb{R}}}\la
h(z-\cdot)\phi',\oz_s\ra^2 dz +\int_0^t\la\sigma\phi^2, \oz_s\ra ds
 \edeqn and
 \begin{eqnarray*}
 & &\int_0^t\int_{S({\mathbb{R}})^{\circ}}
         (e^{-\theta\la\phi,\nu\ra}-1
         +\theta\la\phi,\nu\ra)\hat{N}(ds,d\nu)\cr
 &=&\int_0^tds\int_{{\mathbb{R}}}\oz_s(dx)
    \int_0^{\infty}\gamma(x,d\xi)(e^{-\xi\la\dz_x,\theta\phi\ra}-1
    +\xi\la\dz_x,\theta\phi\ra),\cr
 \end{eqnarray*}
where $\theta>0$ and $\phi\in C^2(\mathbb{R})^{+}$. That is, under
${\bf Q}_{\mu}$ the jump measure $N$ has compensator
 \bgeqn
 \label{f2.35}
 \hat{N}(ds,d\nu)=ds\oz_s(dx)\gamma(x,d\xi)\cdot\dz_{\xi\dz_x}(d\nu)
 ,~~\nu\in
 M({\mathbb{R}}).
 \edeqn
In particular this implies that the jumps of $\omega$ are ${\bf
Q}_{\mu}$-a.s. in $M({\mathbb{R}})$, i.e. positive measures. Observe
that for $\{\phi_i\}_{i=1}^2\subset C^2(\mathbb{R})^+$,
$M_t^c(\phi_1+\phi_2)=M_t^c(\phi_1)+M_t^c(\phi_2)$. According to
($\ref{f2.34}$),
 \begin{eqnarray}
 \label{f2.36}
 \la M^c(\phi_1),M^c(\phi_2)\ra_t
  &=&\frac{1}{2}\int_0^t
  \int_{{\mathbb{R}}^2}\rho(x-y)\phi_1'(x)\phi_2'(y)
     \oz_s(dx)\oz_s(dy)ds\cr
  &&+\frac{1}{2}\int_0^t
  \int_{{\mathbb{R}}^2}\rho(x-y)\phi_2'(x)\phi_1'(y)
     \oz_s(dx)\oz_s(dy)ds\cr
  &&+\int_0^t\la\sigma\phi_1\phi_2,\oz_s\ra ds.
  \end{eqnarray}

\textit{Step 4}. Let
$J_1(\phi,\nu)=\la\phi,\nu\ra1_{\{\la1,\nu\ra\geq
 1\}}$ and $J_2(\phi,\nu)=\la\phi,\nu\ra1_{\{\la1,\nu\ra<
 1\}}$. First, one
 can check that
$$
{\bf Q}_{\mu}\bigg{[}\int_0^{t}\int
 J_1(\phi,\nu)\hat{N}(ds,d\nu)\bigg{]}
 <\infty\quad \textrm{and}\quad{\bf Q}_{\mu}\bigg{[}\int_0^{t}\int
J_2(\phi,\nu)^2\hat{N}(ds,d\nu)\bigg{]}
 <\infty
$$
for $\phi\in C^2(\mbb R)^+$. Then following the argument in Section
2.3 of \cite{[LM05]} we obtain the martingale property of
$M^d(\phi)$.  By Proposition 2.1.28 and Theorem 2.1.33 of
\cite{[JS87]} we can deduce that
\begin{eqnarray*}
\int_0^{t+}\int_{M({\mathbb{R}})^{\circ}}
J_1(\phi,\nu)\tilde{N}(ds,d\nu)
&=&\int_0^{t+}\int_{M({\mathbb{R}})^{\circ}}
J_1(\phi,\nu)N(ds,d\nu)\\&& -\int_0^t\int_{M({\mathbb{R}})^{\circ}}
J_1(\phi,\nu)\hat{N}(ds,d\nu),\quad t\geq0,
\end{eqnarray*}
is a martingale and
$$
\int_0^{t+}\int_{M({\mathbb{R}})^{\circ}}J_2(\phi,\nu)
\tilde{N}(ds,d\nu),\quad t\geq0,
$$
is a square-integrable martingale with quadratic variation process
given by
$$
\la \int_0^{\cdot+}\int_{M({\mathbb{R}})^{\circ}}J_2(\phi,\nu)
\tilde{N}(ds,d\nu)\ra_t= \int_0^{t}\int_{M({\mathbb{R}})^{\circ}}
J_2(\phi,\nu)^2\hat{N}(ds,d\nu).
$$
Recall that
$$
M_t^c(\phi)=M_t(\phi)-M_t^d(\phi).
$$
The fact that  both $M^d(\phi)$ and $M(\phi)$ above are martingales
yields the martingale property of $M^c(\phi)$. We are done. \qed
\begin{lemma}
\label{L2.3} Let ${\bf Q}_{\mu}$ be a probability measure on
$(\Omega,\cal F)$ such that it is a solution of the $({\cal
L},\mu)$-martingale problem. Then
$${\bf Q}_{\mu}[\sup_{0\leq s\leq t}\la1,\oz_s\ra]<\infty.$$
\end{lemma}
{\bf Proof.} According to Theorem \ref{T2.3} and \textit{Step 4} in
its proof, we have
$$
\la1,\oz_t\ra=\la1,\mu\ra+M_t^c(1)+\int_0^t
\int_{M({\mathbb{R}})^{\circ}} \la1,\nu\ra\tilde{N}(ds,d\nu)
$$
is a martingale and we obtain
\begin{eqnarray*}
{\bf Q}_{\mu}\bigg{[}\sup_{0\leq s\leq
t}\la1,\oz_s\ra\bigg{]}\ar\leq\ar\la1,\mu\ra+{\bf
Q}_{\mu}\bigg{[}\sup_{0\leq s\leq t}|M_s^c(1)|\bigg{]}+{\bf
Q}_{\mu}\bigg{[}\sup_{0\leq s\leq t}|\int_0^s\int
J_2(1,\nu)\tilde{N}(ds,d\nu)|\bigg{]}\cr \ar\ar+{\bf
Q}_{\mu}\bigg{[}\sup_{0\leq s\leq t}\int_0^s\int
J_1(1,\nu)N(ds,d\nu)\bigg{]}\cr\ar\ar+{\bf
Q}_{\mu}\bigg{[}\sup_{0\leq s\leq t}\int_0^s\int
J_1(1,\nu)\hat{N}(ds,d\nu)\bigg{]}\cr\ar\leq\ar \la1,\mu\ra+4{\bf
Q}_{\mu}[C_t(1)]+2+{\bf Q}_{\mu}\bigg{[}\sup_{0\leq s\leq
t}[\int_0^s\int J_2(1,\nu)\tilde{N}(ds,d\nu)]^2\bigg{]}\cr
\ar\ar+2\sup_x\int_1^{\infty}\xi\gamma(x,d\xi) \int_0^t{\bf
Q}_{\mu}[\la1,\oz_s\ra] ds\cr \ar\leq\ar
\la1,\mu\ra+2+4||\sigma||\int_0^t{\bf Q}_{\mu}[\la1,\oz_s\ra]
ds+4\sup_x\int_0^1\xi^2 \gamma(x,d\xi)\int_0^t{\bf
Q}_{\mu}[\la1,\oz_s\ra ]ds \cr\ar\ar+
2\sup_x\int_1^{\infty}\xi\gamma(x,d\xi)\int_0^t{\bf
Q}_{\mu}[\la1,\oz_s\ra ]ds\cr\ar\leq\ar
\la1,\mu\ra+2+C_2(\sigma,\gamma)\la1,\mu\ra t,
\end{eqnarray*}
where
$C_2(\sigma,\gamma):=4||\sigma||+2\sup_x\int_1^{\infty}\xi\gamma(x,d\xi)+
4\sup_x\int_0^1\xi^2\gamma(x,d\xi)$ and the second and the third
inequalities follow from  Doob's inequality and the elementary
inequality $|x|\leq x^2+1$. We complete the proof. \qed

\smallskip \noindent In accordance with the notation used in Theorem
$\ref{T2.3}$, set
$$X^L_t:=\oz_t-\int_0^{t+}\int_{M(\mbb
R)^{\circ}}\nu\cdot1_{\{\la1,\nu\ra\geq l\}}N(ds,d\nu).
$$
By Theorem $\ref{T2.3}$,
\begin{eqnarray}
\label{f2.37}
 \la \phi,X^L_t\ra=\la \phi,\mu\ra+\int_0^t\la
 a\phi'',\oz_s\ra ds+M_t^c(\phi)\ar+\ar\int_0^{t+}\int_{M(\mbb
R)^{\circ}}\la\phi,\nu\ra1_{\{\la1,\nu\ra<l\}}\tilde{N}(ds,d\nu)\cr
\ar-\ar\int_0^{t+}\int_{M(\mbb
R)^{\circ}}\la\phi,\nu\ra1_{\{\la1,\nu\ra\geq l\}}\hat{N}(ds,d\nu).
\end{eqnarray}
 Thus if
$F(\mu)=f(\la\phi_1,\mu\ra,\cdots,\la\phi_m,\mu\ra)\in \cal D(\cal
L)$, then by It\^{o}'s formula
 \begin{eqnarray*}
I_t:= F(X^L_t)\ar+\ar\frac{1}{2}\sum_{i=1}^{m}\int_0^tds\int_{\mbb
R}\oz_s(dx)\int_l^{\infty}\gamma(x,d\xi)
f^i(\la\phi_1,X^L_s\ra,\cdots,\la\phi_n,X^L_s\ra)\xi\phi_i(x)\cr
    \ar-\ar\frac{1}{2}\sum_{i=1}^{m}
    \int_0^tf^i(\la\phi_1,X^L_s\ra,\cdots,\la\phi_n,X^L_s\ra)\la
    a\phi_i'',\oz_s\ra ds\cr
 \ar-\ar\frac{1}{2}\sum_{i,j=1}^m\int_0^t
 f^{ij}(\la\phi_1,X^L_s\ra,\cdots,\la\phi_n,X^L_s\ra)d\la
 M^c(\phi_i),M^c(\phi_j)\ra s\cr
 \ar-\ar\int_0^tds\int_{{\mathbb{R}}}\oz_s(dx)
 \int_0^{l}\gamma(x,d\xi)
 \big{\{}f(\la\phi_1,X^L_s\ra+\xi\phi_1(x),\cdots,\la\phi_n,X^L_s\ra
 +\xi\phi_n(x))\cr
 &  &~~~~~~-f(\la\phi_1,X^L_s\ra,\cdots,\la\phi_n,X^L_s\ra)
 -\xi\sum_{i=1}^m\phi_i(x)
 f^i(\la\phi_1,X^L_s\ra,\cdots,\la\phi_n,X^L_s\ra)\big{\}}
\end{eqnarray*}
is a local martingale under $\bf Q _{\mu}$.

\noindent Let $\tau^1=\inf\{t\geq0: \la1,\oz_t\ra\geq
l+\la1,\mu\ra\}\wz T$ and $\tau^2=\inf\{t\geq0:
|\la1,\oz_t\ra-\la1,\oz_{t-}\ra|\geq l\}$. Set
$\tau=\tau^1\wz\tau^2$. The following lemma gives  another
martingale characterization for $X^L$.

\begin{lemma}
\label{L2.4} Let $\bf P_{\mu}$ be a  probability measure on
$(\Omega,\cal F)$ such that ${\bf P}  _{\mu}(\omega_0=\mu)=1$. Then
 \begin{eqnarray} \label{f2.38} I_t(\phi):=\exp\bigg{\{}-\la
\phi,X^L_{t\wz\tau}\ra\ar+\ar\int_0^{t\wz\tau}[\la
a\phi'',\oz_s\ra-\int_{\mbb R}\la h(z-\cdot)\phi',\oz_s\ra^2dz]ds\cr
\ar-\ar\int_0^{t\wz\tau}ds\int_{\mbb
R}\oz_s(dx)\int_{l}^{\infty}\xi\phi(x)\gamma(x,d\xi)\cr
\ar-\ar\int_0^{t\wz\tau}ds\int_{\mbb
R}\oz_s(dx)\int_{0}^{l}(e^{-\xi\phi(x)}-1+\xi\phi(x))\gamma(x,d\xi)
\bigg{\}}
\end{eqnarray}
is a $\bf P_{\mu}$-martingale for every $\phi\in C^2(\mbb R)^{++}$
if and only if  $\{I_{t\wz\tau}\}$ is a $\bf P_{\mu}$-martingale for
each $F\in \cal D(\cal L)$.
\end{lemma}
$\bf Proof$.   The desired result follows from the formula of
integration by parts and  the same argument as in the proof of
Th\'{e}or\`{e}m 7 of \cite{[ER91]}. \qed

\noindent The next two theorems are analogous to  Theorem (3.1) and
Theorem (3.3) of \cite{[S75]}.
\begin{theorem}
\label{T2.4} Given a probability measure $\bf P$ on $(\Omega,\cal
F)$ such that ${\bf P}(\omega(0)=\mu)=1$ and $\{I(t\wz\tau):t\geq
0\}$ is a ${\bf P}$-martingale. Define
$$
{\bf S}_{\omega}=\dz_{\omega}\otimes {\bf Q}'_ {\{\tau(\omega),
X^L_{\tau(\omega)}\}}
$$
and
$$
{\bf P}'(A)={\bf P}[{\bf S}_{\omega}(A)],~~A\in \cal F,
$$
where ${\bf S}_{\omega}$ is a measure on $(\Omega,\cal F)$
satisfying
$$
{\bf S}_{\omega}(A_1\cap A_2)=1_{A_1}(\omega){\bf Q}'_
{\{\tau(\omega), X^L_{\tau(\omega)}\}}(A_2)
$$
for $A_1\in \sigma(\bigcup_{0\leq s<\tau(\omega)}{\cal F}_s)$ and
$A_2\in\cal F^{\tau(\omega)}$. Define ${\cal
F}_{\tau-}=\sigma\{X^L_{t\wz\tau}: t\geq0\}$. Then ${\bf P}'$ is
also a solution of $({\cal L}',\mu)$-martingale problem and ${\bf
P}={\bf Q}_{\mu}'$ on ${\cal F}_{\tau-}$. In particular, we can take
$\bf P={\bf Q}_{\mu}.$
\end{theorem}
\noindent${\bf Proof.}$ Let $0\leq t_1<t_2$ and $A\in{\cal
F}_{t_1}$. Given $\oz\in\Omega$, for this proof only, let $y(t,\oz)$
denote the position of $\oz$ at time $t$ for convenient. Let $F\in
{\cal D}(\cal L)$. Then
\begin{eqnarray*}
{\bf P}'[1_AF(y_{t_2})]\ar=\ar{\bf
P}[1_{A\cap\{\tau>t_2\}}F(X^L_{t_2})]+{\bf P}[1_{A\cap\{t_1<\tau\leq
t_2\}}{\bf Q}'_{\tau(\oz),X^L_{\tau(\oz)}}[F(y_{t_2})]]\cr \ar\ar+
{\bf P}[1_{\{\tau\leq t_1\}}{\bf
S}_{\oz}[1_AF(y_{t_2})]]=I_1+I_2+I_3.
\end{eqnarray*}
By the martingale formula of $\bf Q'$
$$
I_2={\bf P}[1_{A\cap\{t_1<\tau\leq t_2\}}F(X^L_{\tau})]+{\bf
P}'[1_{A\cap\{t_1<\tau\leq t_2\}}\int_{\tau}^{t_2}{\cal
L}'F(y_u)du],
$$
and
\begin{eqnarray*}
I_1+I_2\ar=\ar{\bf P}[1_{A\cap\{\tau>t_1\}}F(X^L_{\tau\wz
t_2})]+{\bf P}'[1_{A\cap\{t_1<\tau\leq t_2\}}\int_{\tau}^{t_2}{\cal
L'} F(y_u)du]\cr \ar=\ar{\bf
P}[1_{A\cap\{\tau>t_1\}}F(X^L_{t_1})]+{\bf
P}[1_{A\cap\{\tau>t_1\}}\int_{t_1}^{\tau\wz t_2}{\cal
L}'F(y_u)du]\cr\ar\ar+{\bf P}'[1_{A\cap\{\tau>t_1\}}\int_{\tau\wz
t_2}^{t_2}{\cal L}'F(y_u)du]\cr\ar=\ar{\bf
P}'[1_{A\cap\{\tau>t_1\}}F(y_{t_1})]+{\bf
P}'[1_{A\cap\{\tau>t_1\}}\int_{t_1}^{t_2}{\cal L}'F(y_u)du],
\end{eqnarray*}
where the second equality follows from that $\{I_{t\wz\tau}\}$ is a
martingale and the fact that $F(X^L_t)-I_t=\int_0^t {\cal L'}
F(\omega_s)ds$ for $ \tau>t$. On the other hand,
$$
I_3={\bf P}'[1_{A\cap\{\tau\leq t_1\}}F(y_{t_1})]+{\bf
P}'[1_{A\cap\{\tau\leq t_1\}}\int_{t_1}^{t_2}{\cal L}'F(y_u)du].
$$
Thus  ${\bf P}'$ solves the $({\cal L}', \mu)$-martingale problem.
Then the desired conclusion follows from the uniqueness of the
$({\cal L}', \mu)$-martingale problem. \qed

\begin{theorem}
\label{T2.5}Let $M_l(\mbb R)=\{\nu:\la1,\nu\ra\geq l\}$. There is a
${\cal F}_{\tau-}$-measurable function
$\tau':\Omega\rightarrow[0,T]$ such that for $\Gamma\in
\mathfrak{B}(M_l(\mbb R))$, \beqlb
 \label{f2.39}
 {\bf Q}_{\mu}[\int_0^{\tau+}N(ds,\Gamma)|{\cal
 F}_{\tau-}]=\int_0^{\tau'}   \exp\{-\int_0^tds\int_{\mbb
R}X^L_{s\wz\tau}(dx)\int_l^{\infty}\gamma(x,d\xi)\}
K(X^L_{t\wz\tau},\Gamma) dt
 \eeqlb
  holds for any solution
${\bf Q}_{\mu}$ to the $(\cal L,\mu)$-martingale problem. In
particular, ${\bf Q}_{\mu}$ is uniquely determined on ${\cal
F}_{\tau}$.
\end{theorem}
{\bf Proof.} In accordance with the notation used in Theorem
$\ref{T2.3}$, we have
 \bgeqn
 \label{f2.40}
 \int_0^{t+}N(ds,\Gamma)=\int_0^{t+}\tilde{N}(ds,\Gamma)+
 \int_0^{t}\hat{N}(ds,\Gamma),
 \edeqn
 where $\hat{N}(ds,\Gamma)$ is determined by (\ref{f2.35}).
An application of It\^{o}'s formula and integration by parts shows
that
$$J_t^{\alpha}:=\exp[\alpha\int_0^{t+}N(ds,\Gamma)-
\int_0^{t}(e^{\alpha}-1)\hat{N}(ds,\Gamma)]$$ is a ${\bf
Q}_{\mu}$-martingale for all $\alpha\in \mbb R$. Combing
(\ref{f2.37}) and (\ref{f2.40}) together and using It\^{o}'s formula
and integration by parts again we see $I_t(\phi)J_t^{\alpha}$ is a
${\bf Q}_{\mu}$-martingale for all $\phi\in C^2(\mbb R)^{++}$. By
Theorem \ref{T2.4} and Lemma \ref{L2.4}, $I_t(\phi)$,
$J_t^{\alpha}$, $\bf Q'$, ${\bf Q}_{\mu}$ and ${\cal F}_{\tau-}$
satisfy the requirement of Theorem (3.2) in \cite{[S75]}. Hence, for
any bounded stopping time ${t_0}$,
 \bgeqn \label{f2.41}{\bf Q}_{\mu}[J_{{t_0}}^{\alpha}|{\cal
F}_{\tau-}]=1 ~~~(a.s.,{\bf Q}_{\mu}).
 \edeqn
 Since $\tau^1$ is a stopping time and $\tau^1\leq T$, we can find a
 measurable function $f:(M(\mbb R))^{{\bf N}}\rightarrow[0,T]$ and $0\leq
 t_1<\cdots<t_n<
 \cdots\leq T$ such that
 $$
 \tau^1=f(\oz_{t_1},\cdots,\oz_{t_n},\cdots).
 $$
Define
$$
\tau'=f(X^L_{t_1\wz\tau},\cdots,X^L_{t_n\wz\tau},\cdots).
$$
Note that $\tau^1=\tau'$ if $\tau^1<\tau^2$.  On the other hand,
\begin{eqnarray*}
{\bf Q}_\mu[\tau\leq t|{\cal F}_{\tau-}]\ar=\ar1_{[0,t]}(\tau'){\bf
Q}_{\mu}[\tau^2>\tau^1|{\cal F}_{\tau-}]+{\bf
Q}_{\mu}[\tau^2\leq\tau^1\wz t|{\cal
F}_{\tau-}]\cr\ar=\ar1_{[0,t]}(\tau'){\bf
Q}_{\mu}[1-\int_0^{\tau+}N(ds,M_l(\mbb R))|{\cal
F}_{\tau-}]\cr\ar\ar+{\bf Q}_{\mu}[\int_0^{(t\wz\tau)+}N(ds,M_l(\mbb
R))|{\cal F}_{\tau-}].
\end{eqnarray*}
According to (\ref{f2.41}),
 \beqlb
 \label{f2.42}
 \ar\ar{\bf Q}_{\mu}[\int_0^{(t\wz\tau)+}N(ds,\Gamma)|{\cal
 F}_{\tau-}]\cr\ar=\ar{\bf
 Q}_{\mu}[\int_0^{t\wz\tau}\hat{N}(ds,\Gamma)|{\cal
 F}_{\tau-}]\cr\ar=\ar\int_0^t{\bf Q}_{\mu}[\tau>s|{\cal
 F}_{\tau-}]\int_{\mbb
 R}X^L_{s\wz\tau}(dx)\int_0^{\infty}\gamma(x,d\xi)1_{\{\xi\dz_x\in
 \Gamma\}}ds
 \eeqlb
for any $\Gamma\in\mathfrak{B}(M_l(\mbb R))$. Thus
 \beqnn
 {\bf Q}_{\mu}[\tau\leq t|{\cal
 F}_{\tau-}]\ar=\ar1_{[0,t]}(\tau')\left(1-\int_0^t{\bf Q}_{\mu}
 [\tau>s|{\cal
 F}_{\tau-}]\int_{\mbb
 R}X^L_{s\wz\tau}(dx)\int_0^{\infty}\gamma(x,d\xi)1_{\{\xi\dz_x\in
 M_l(\mbb R)\}}ds\right)\cr\ar\ar+\int_0^t{\bf Q}_{\mu}[\tau>s|{\cal
 F}_{\tau-}]\int_{\mbb
 R}X^L_{s\wz\tau}(dx)\int_0^{\infty}\gamma(x,d\xi)1_{\{\xi\dz_x
 \in M_l(\mbb R)\}}ds
 \eeqnn
and so
$$
{\bf Q}_{\mu}[\tau>t|{\cal
F}_{\tau-}]=1_{(t,\infty)}(\tau')\exp\{-\int_0^tds\int_{\mbb
R}X^L_{s\wz\tau}(dx)\int_l^{\infty}\gamma(x,d\xi)\}.
$$
Plugging this back into (\ref{f2.42}) and setting $t=T$, we obtain
(\ref{f2.39}).

\noindent Finally, since
$\oz_{\tau}=X^L_{\tau}+\int_0^{\tau+}\int\nu1_{\{\la1,\nu\ra\geq
l\}}N(ds,d\nu)$, we see that  the distribution of $\oz_{\tau}$ under
${\bf Q}_{\mu}$ given ${\cal F}_{\tau-}$ is uniquely determined,
and, therefore ${\bf Q}_{\mu}$ is uniquely determined on ${\cal
F}_{\tau}$.  \qed

\begin{lemma}
\label{L2.5} Let ${\bf Q}_{\mu}$ be a solution of $({\cal
L},\mu)$-martingale problem. Given a finite stopping time $\beta$,
let ${\cal Q}_{\omega}$ be a regular conditional probability
distribution of ${\bf Q}_{\mu}|{\cal F}_\beta$. Then there is an
$N\in {\cal F}_{\beta}$ such that ${\bf Q}_{\mu}(N)=0$ and when
$\oz\notin N$
$$
F(\oz'_{t\vee\beta(\oz)})-F(\oz'_{\beta(\oz)})
-\int_{\beta(\oz)}^{t\vee\beta(\oz)}{\cal
L}F(\oz'_s)ds
$$
under  ${\cal Q}_{\oz}$  is a martingale for $F\in {\cal D}_0(\cal
L)$. In particular, it is a local martingale for all $F\in\cal
D(\cal L)$.
\end{lemma}
$\bf Proof.$ The argument in this proof is exactly the same as that
in Theorem 6.1.3 of \cite{[SV79]}. We omit it here. \qed

\noindent Now, we come to our main theorem in this subsection.
\begin{theorem}
\label{T2.6} Suppose that for $l>1$, the $({\cal
L}',\mu)$-martingale problem is well-posed. Then uniqueness hold for
$({\cal L},\mu)$-martingale problem.
\end{theorem}
{\bf Proof}. Suppose ${\bf Q}_{\mu}$ is a solution of $({\cal
L},\mu)$-martingale problem. Define $\beta_0=0$ and
$$
\beta_{n+1}=\left(\inf\{t\geq \beta_n:|\la1,
\oz_t\ra-\la1,\oz_{t-}\ra|\geq l\textrm{ or
}\la1,\oz_t\ra-\la1,\oz_{\beta_n}\ra\geq l\}\right)\wedge
(\beta_n+l).
$$
 Then for each $n\geq1$, $\beta_n$ is a stopping time bounded by $nl$.
 By Lemma \ref{L2.5} and Theorem \ref{T2.5}, we can prove by
induction that ${\bf Q}_{\mu}$ is uniquely determined on ${\cal
F}_{\beta_n}$ for all $n\geq 1$. In order to get the desired
conclusion we only need to show that ${\bf Q}_{\mu}(\beta_n\leq
t)\rightarrow0$ as $n\rightarrow\infty$ for each $t>0$.

\noindent Let $\beta^1_0=0$ and $\beta^2_0=0$. Define
$$\beta^1_{n+1}=\inf\{t\geq
\beta^1_n:\la1,\oz_t\ra-\la1,\oz_{\beta^1_n}\ra\geq l\}
$$ and
$$\beta^2_{n+1}=\inf\{t\geq \beta^2_n:\la1,\oz_t\ra-\la1,\oz_{t-}\ra\geq
l\}.
$$
It is easy to see that in order to get the desired conclusion it
suffices to show that  ${\bf Q}_{\mu}(\beta_n^1\leq t)\rightarrow0$
and ${\bf Q}_{\mu}(\beta_n^2\leq t)\rightarrow0$ as
$n\rightarrow\infty$. First, by Lemma \ref{L2.3}, we can deduce that
$$
\lim_{n\rightarrow\infty}{\bf Q}_{\mu}(\beta_n^1\leq t)=0.
$$
Then
\begin{eqnarray*} \sum_{0<s\leq t}1_{\{\la1,\Delta\oz_s\ra\geq
l\}}\ar\leq\ar \sum_{0<s\leq
t}\la1,\Delta\oz_s\ra1_{\{\la1,\Delta\oz_s\ra\geq l\}}\cr
 \ar=\ar\int_0^{t+}\int
 _{M(\mbb R)^{\circ}}\la1,\nu\ra1_{\{\la1,\nu\ra\geq l\}}N(ds,d\nu).
\end{eqnarray*}
But according to the \textit{Step 4} in the proof of Theorem
\ref{T2.3},
$${\bf Q}_{\mu}\bigg{[}\int_0^{t+}\int
 _{M(\mbb R)^{\circ}}\la1,\nu\ra1_{\{\la1,\nu\ra\geq
 l\}}N(ds,d\nu)\bigg{]}<\infty,$$
which yields that
 $$\lim_{n\rightarrow\infty}{\bf
Q}_{\mu}(\beta_n^2\leq t)=0.$$
 \qed

\section{Existence}
\subsection{Interacting-branching particle system }
We first give a formulation of the interacting-branching particle
system. Then we construct a solution of the $(\cal
L',\mu)$-martingale problem by using particle system approximation.
We recall that
$$
G^m:=\frac{1}{2}\sum_{i=i}^m a(x_i)\frac{\partial^2}{\partial
x_i^2} +\frac{1}{2}\sum_{i,j=1,i\neq
j}^m\rho(x_i-x_j)\frac{\partial^2}{\partial x_i\partial x_j}.
$$
Suppose that $X_t=(x_1(t),\cdots,x_m(t))$ is a Markov process in
$\mathbb{R}^m$ generated by ${G}^m$. By Lemma 2.3.2 of \cite{[D93]}
we know that $X_t=(x_1(t),\cdots,x_m(t))$ is an exchangeable Feller
process.  Let $N(\mathbb{R})$ denote the space of integer-valued
measures on $\mathbb{R}$. For $\theta>0$, let
$M_{\theta}(\mathbb{R})=\{\theta^{-1}\sigma:\sigma\in
N(\mathbb{R})\}.$ Let $\zeta$ be the mapping from
$\cup
_{m=1}^{\infty}\mathbb{R}^m$ to $M_{\theta}(\mathbb{R})$
defined by
$$
\zeta(x_1,\cdots,x_m)=\frac{1}{\theta}\sum_{i=1}^m\delta_{x_i},
\textrm{\ \ \ }m\geq 1.
$$
By Proposition 2.3.3 of \cite{[D93]} we know that $\zeta(X_t)$ is a
Feller Markov process in $M_{\theta}(\mathbb{R})$ with generator
$\mathcal {A}_{\theta}$ given by
 \begin{eqnarray}
 \label{3.1}
 \mathcal {A}_{\theta}F(\mu)&=&\frac{1}{2}\int_{\mathbb{R}}a(x)
 \frac{d^2}{dx^2}\frac{\delta F(\mu)}{\delta\mu(x)}\mu(dx)
 +\frac{1}{2\theta}\int_{\mathbb{R}^2}c(x)c(y)\frac{d^2}{dx dy}
 \frac{\delta^2F(\mu)}{\delta\mu(x)\delta\mu(y)}\delta_x(dy)\mu(dx)\cr
 &&+\frac{1}{2}\int_{\mathbb{R}^2}\rho(x-y)\frac{d^2}{dx dy}
 \frac{\delta^2F(\mu)}{\delta\mu(x)\delta\mu(y)}\mu(dx)\mu(dy).
 \end{eqnarray}
In particular, if
 \begin{eqnarray}\label{3.1FormP}F(\mu)=f\left(\la\phi_1,\mu\ra,
 \cdots,\la\phi_n,\mu\ra\right),\textrm{\ \ \
 }\mu\in M_\theta(\mathbb{R}),\end{eqnarray}
 for $f\in C^2(\mathbb{R}^n)$ and $\{\phi_i\}\subset
 C^2(\mathbb{R})$, then
\begin{eqnarray}
\label{3.2} \mathcal {A}_\theta F(\mu) &=&\frac{1}{2}\sum_{i=1}^nf^i
\left(\la\phi_1,\mu\ra,\cdots,\la\phi_n,\mu\ra\right)\la
         a\phi''_i,\mu\ra \cr
&&+\frac{1}{2\theta}\sum_{i,j=1}^n
   f^{ij}\left(\la\phi_1,\mu\ra,\cdots,\la\phi_n,\mu\ra\right)\la
    c^2\phi'_i\phi'_j,\mu\ra\cr
&&+\frac{1}{2}\sum_{i,j=1}^nf^{ij}
\left(\la\phi_1,\mu\ra,\cdots,\la\phi_n,\mu\ra\right)
   \int_{\mathbb{R}^2}\rho(x-y)\phi'_i(x)\phi'_j(y)\mu(dx)\mu(dy).
\end{eqnarray}

\noindent Now we introduce a branching mechanism to the interacting
particle system. Suppose that for each $x\in\mathbb{R}$ we have a
discrete probability distribution $p(x)=\{p_i(x):i=0,1,\cdots\}$
such that each $p_i(\cdot)$ is a Borel measurable function on
$\mathbb{R}$. This serves as the distribution of the offspring
number produced by a particle that dies at site $x\in \mathbb{R}$.
We assume that
 \bgeqn
 \label{3.3}
 \sum_{i=1}^{\infty}ip_i(x)\leq1,
 \edeqn
and
 \bgeqn
 \label{3.4}
 \sigma_p(x):=\sum_{i=1}^{\infty}i^2 p_i(x)
 \edeqn
is bounded in $x\in\mathbb{R}$. For $0\leq z\leq1$, let
 \bgeqn
 \label{3.5}
 g(x,z):=\sum_{i=0}^{\infty}p_i(x)z^i.
 \edeqn
Let $\Gamma_{\theta}(\mu,d\nu)$ be the probability kernel on
$M_\theta(\mathbb{R})$ defined by
 \bgeqn
 \label{3.6}
 \int_{M_{\theta}(\mathbb{R})}F(\nu)\Gamma_{\theta}(\mu,d\nu)
 =\frac{1}{\la1,\mu\ra}\left\la\sum_{j=0}^{\infty}
 p_j(x)F\left(\mu+(j-1)\theta^{-1}\delta_{x}\right),\mu\right\ra,
 \edeqn
where $\mu\in M_{\theta}(\mathbb{R})$ is given by
$$ \mu=\frac{1}{\theta}\sum_{i=1}^{\theta\la1,\mu\ra}\delta_{x_i}.$$
For a constant $\lz>0$, we define the bounded operator
$\mathcal{B}_{\theta}$ on $B(M_{\theta}(\mathbb{R})$ by
 \bgeqn
 \label{3.7}
 \mathcal{B}_{\theta}F(\mu)=
 \lambda\theta(\theta\wedge\la1,\mu\ra)\int_{M_{\theta}(\mathbb{R})}
 [F(\nu)-F(\mu)]\Gamma_{\theta}(\mu,d\nu).
 \edeqn
For $\mathcal{A}_{\theta}$ generates a  Markov process on
$M_{\theta}(\mathbb{R})$, then
$\mathcal{L}_{\theta}:=\mathcal{A}_{\theta}+\mathcal{B}_{\theta}$
also generates a Markov process; see Problem 4.11.3 of
\cite{[EK86]}. By martingale inequality and Theorem 4.3.6 of
\cite{[EK86]}, we obtain  that the corresponding Markov process has
a modification with sample paths in $D([0,\infty),
M_{\theta}(\mathbb{R}))$. We shall call the process generated by
$\mathcal{L}_{\theta}$ an interacting-branching particle system with
parameter $(a, \rho, \gamma, \lambda,p)$ and unit mass $1/\theta$.

\subsection{Particle system approximation}
 \noindent Recall that
 \bgeqn
 \label{f3.9}
\Psi_0(x,z):=\frac{1}{2}\sigma(x)z^2+\int_l^{\infty}\xi\gamma(x,d\xi)z
  +\int_0^l(e^{-z\xi}-1+z\xi)\gamma(x,d\xi).
 \edeqn
According to the conditions (i) and (iii) on the $\sigma$ and
$\gamma(x,d\xi)$, $\Psi_0(x,\phi(x))\in C(\mathbb{R})$ can be
extended continuously to $\hat{\mathbb{R}}$ for $\phi\in
C^2_{\partial}(\mathbb{R})^{++}$. And, if
 \bgeqn
 \label{f3.10}
 F(\mu)=f(\la\phi_1,\mu\ra,\cdots,\la\phi_n,\mu\ra),~~~~\mu\in
 M(\mathbb{R}),
 \edeqn
for $f\in C^2(\mathbb{R}^n)$ and $\{\phi_i\}\subset
C^2(\mathbb{R})$, then
 \begin{eqnarray}
 \label{f3.11}
 \mathcal{A}F(\mu)&=&\frac{1}{2}\sum_{j=1}^n
 f^i(\la\phi_1,\mu\ra,\cdots,\la\phi_n,\mu\ra)\la
 a\phi''_i,\mu\ra\cr
 &&+\frac{1}{2}\sum_{i,j=1}^n
 f^{ij}(\la\phi_1,\mu\ra,\cdots,\la\phi_n,\mu\ra)
 \int_{\mathbb{R}^2}\rho(x-y)\phi_i'(x)\phi_j'(y)\mu^2(dxdy)
 \end{eqnarray}
and
 \begin{eqnarray}
 \label{f3.12}
 \mathcal{B}'F(\mu)&=&\frac{1}{2}\sum_{i,j=1}^n
 f^{ij}(\la\phi_1,\mu\ra,\cdots,\la\phi_n,\mu\ra)
 \la\sigma\phi_i\phi_j,\mu\ra\cr
 &&-\int_{\mbb R}\mu(dx)\int_l^{\infty}\xi\gamma(x,d\xi)\sum_{i=1}^n
  f^i(\la\phi_1,\mu\ra,\cdots,\la\phi_n,\mu\ra)\phi_i(x)\cr
 &&+\int_{\mathbb{R}}\mu(dx)\int_0^{l}\{
 f(\la\phi_1,\mu\ra+\xi\phi_1(x),\cdots,
 \la\phi_n,\mu\ra+\xi\phi_n(x))\cr
 &&~~~-f(\la\phi_1,\mu\ra,\cdots,\la\phi_n,\mu\ra)
  -\xi\sum_{i=1}^n
  f^i(\la\phi_1,\mu\ra,\cdots,\la\phi_n,\mu\ra)\phi_i(x)\}
  \gamma(x,d\xi).
 \end{eqnarray}
Suppose $\{X_t^{(k)} : t\geq0\}$ is a sequence of
$c\acute{a}dl\acute{a}g$ interacting-branching particle systems with
parameters $(a,\rho, \gamma, \lambda_k, p^{(k)})$ and unit mass
$1/k$ and initial states $X_0^k=\mu_k\in M_k(\mathbb{R})$. We can
regard $\{X_t^{(k)} : t\geq 0\}$ as a process with state space
$M(\hat{\mathbb{R}})$. Let $\sigma_p^k$ and $g_k$ be defined by
($\ref{3.4}$) and ($\ref{3.5}$) respectively with $p_i$ replaced by
$p_i^{(k)}$. Let
 \bgeqn
 \label{f3.13}
 \psi_k(x,z):=k\lz_k[g_k(x, 1-z/k)-(1-z/k)],
 ~~0\leq z\leq k.
 \edeqn
We have that
$\frac{d}{dz}\psi_k(x,0+)=\lz_k[1-\frac{d}{dz}g_k(x,1)]$ and
$\frac{d^2}{dz^2}\psi_k(x,0+)=\lz_k\sigma_p^k/k$.
\smallskip

\begin{lemma} \label{L3.1}Suppose that the sequence
$\{\lz_k\sigma_p^k/k\}$ and $\{\la1,\mu_k\ra\}$ are bounded. Then
$\{X_t^{(k)} : t\geq0 \}$ form a tight sequence in $D([0,+\infty),
M(\hat{\mathbb{R}}))$.\end{lemma}
 \textbf{Proof}. By ($\ref{3.3}$),
it is easy to see that $\{\la1, X_t^{(k)}\ra:t\geq0\}$ is a
supermartingale. By using martingale inequality, one can check that
$\{X_t^{(k)}:t\geq0\}$ satisfies the compact containment condition.
Let $\mathcal{L}_k$ denote the generator of $\{X_t^{(k)}:t\geq0\}$
and let $F$ be given by ($\ref{f3.10}$) with $f\in
C_0^2(\mathbb{R}^n)$ and with each $\phi_i\in
C_{\partial}^2(\mathbb{R})^{++}$. Then
$$
F(X_t^{(k)})-F(X_0^{(k)})-\int_0^t\mathcal{L}_kF(X_s^{(k)})ds,~~t\geq0,
$$
is a martingale and the desired tightness result follows from
Theorem 3.9.4
 of  Ethier and Kurtz \cite{[EK86]}.\qed

\noindent In the sequel of this subsection, we assume
$\{\phi_i\}\subset C^2_{\partial}(\mathbb{R})$. In this case,
($\ref{f3.10}$), ($\ref{f3.11}$) and ($\ref{f3.12}$) can be extended
to continuous functions on $M(\hat{\mathbb{R}})$. Let
$\hat{\mathcal{A}}F(\mu)$ and $\hat{\mathcal{B}}'F(\mu)$ be defined
respectively by the right hand side of the ($\ref{f3.11}$) and
($\ref{f3.12}$) and let
$\hat{\mathcal{L}}'F(\mu)=\hat{\mathcal{A}}F(\mu)+\hat{\mathcal{B}}'F(\mu),$
all defined as continuous functions on $M(\hat{\mathbb{R}})$.

\begin{lemma}
\label{L3.2} Let $\mathcal{D}_0(\hat{\mathcal{L}}')$ be the totality
of all functions of the form ($\ref{f3.10}$) with $f\in
C_0^2(\mathbb{R}^n)$ and with each $\phi_i\in
C_{\partial}^2(\mathbb{R})^{++}$. Suppose that
$\mu_k\rightarrow\mu\in M(\hat{\mathbb{R}})$  as $k\rightarrow
+\infty$ and the sequence $\{\lz_k\sigma_p^k/k\}$  is bounded. If
for each $h\geq0$, $\psi_k(x,z)\rightarrow\Psi_0(x,z)$ uniformly on
$\mathbb{R}\times[0,h]$ and $\frac{d}{dz}\psi_k(x,0+)\rightarrow
\frac{d}{dz}\Psi_0(x,0)$ uniformly on $\mbb R$ as
$k\rightarrow+\infty$, then for each
$F\in\mathcal{D}_0(\hat{\mathcal{L}}')$,
 \bgeqn
 \label{f3.14}
 F(\omega_t)-F(\omega_0)-\int_0^t\hat{\mathcal{L}}'F(\omega_s)ds,~~~t\geq0,
 \edeqn
is a martingale under any limit point $\mathbf{Q}_{\mu}$ of the
distributions of $\{X_t^{(k)}:t\geq0\}$, where $\{\omega_t:t\geq0\}$
denotes the coordinate process of
$D([0,\infty),M(\hat{\mathbb{R}}))$.
\end{lemma}
\textbf{Proof}. By passing to a subsequence if it is necessary, we
may assume that the distribution of $\{X_t^{(k)}:t\geq0\}$ on
$D([0,+\infty),M(\hat{\mathbb{R}}))$ converges to $\mathbf{Q}_\mu$.
Using Skorokhod's representation, we may assume that the processes
$\{X_t^{(k)}:t\geq0\}$ are defined on the same probability space and
the sequence converges almost surely to a c\`{a}dl\`{a}g process
$\{X_t:t\geq0\}$ with distribution  $\mathbf{Q}_{\mu}$ on
$D([0,\infty),M(\hat{\mathbb{R}}))$ (\cite{[EK86]}, p.102). Let
$K(X)=\{t\geq0:\mathbf{P}\{X_t=X_{t-}\}=1\}$. By Lemma 3.7.7 of
\cite{[EK86]}, the complement of the set $K(X)$ is at most countable
and by Proposition 3.5.2 of \cite{[EK86]}, for each $t\in K(X)$ we
have a.s. $\lim_{k\rightarrow\infty}X_t^{(k)}=X_t$. Our proof will
be divided into 3 steps.

  \emph{Step 1}. We shall
show that
 \bgeqn
 \label{f3.15}
 M_t(\phi):=\la\phi,X_t\ra-\la\phi, X_0\ra-\frac{1}{2}\int_0^t\la
 a\phi'',X_s\ra ds+
 \int_0^tds\int_{\mbb R}X_s(dx)\phi(x)\int_l^{\infty}\xi\gamma(x,d\xi),~~t\geq0,
 \edeqn
 is a  square-integrable martingale with
 $\phi\in C^2_{\partial}(\mathbb{R})$.
 First, Fatou's Lemma
tells us $\mathbf{E}\la
1,X_t\ra\leq\liminf\limits_{k\rightarrow\infty}\mathbf{E}\la
 1,X_t^{(k)}\ra.$
 On the other hand, for $\mu_k\in M_k(\mathbb{R})$ we can get that
 $$\mathcal{L}_k\la\phi,\mu_k\ra=\frac{1}{2}\la a\phi'',\mu_k\ra
 -\frac{k\wz\mu_k(1)}{\mu_k(1)}\la\frac{d}{dz}\psi_k(x,0+)\phi(x),\mu_k\ra.$$
  Then for $t\in K(X)$
  \bgeqn
  \label{f3.16}\mathbf{E}\la1,X_t\ra\leq\liminf_{k\rightarrow\infty}
  \mathbf{E}\la
 1,X_t^{(k)}\ra\leq\liminf_{k\rightarrow\infty}{\bf E}
 \la1,X_0^{(k)}\ra\leq\la1,X_0\ra
  \edeqn
 and a.s.
 $$\lim_{k\rightarrow\infty}\mathcal{L}_k\la\phi,X_t^{(k)}\ra
 =\hat{\mathcal{L}}'\la\phi,X_t\ra=\frac{1}{2}\la a\phi'',X_t\ra
 -\int_{\mbb R}X_t(dx)\phi(x)\int_l^{\infty}\xi\gamma(x,d\xi).$$
 Suppose that $\{H_i\}_{i=1}^n\subset C(M(\hat{\mathbb{R}}))$ and
 $\{t_i\}_{i=1}^{n+1}\subset K(X)$ with $0\leq
 t_1<\cdots<t_n<t_{n+1}$. Then
  \begin{eqnarray*}
  &&\mathbf{E}\big{\{}\big{[}\la\phi,X_{t_{n+1}}\ra-\la\phi,X_{t_n}\ra
  -\int_{t_n}^{t_{n+1}}
    \hat{\mathcal{L}}'\la\phi,X_s\ra ds\big{]}\prod_{i=1}^nH_i(X_{t_i})\big{\}}\cr
  &=&\mathbf{E}\big{\{}\la\phi,X_{t_{n+1}}\ra\prod_{i=1}^nH_i(X_{t_i})\big{\}}-
     \mathbf{E}\big{\{}\la\phi,X_{t_n}\ra\prod_{i=1}^nH_i(X_{t_i})\big{\}}\cr
  & &~~   -\int_{t_n}^{t_{n+1}}
      \mathbf{E}\big{\{}\hat{\mathcal{L}}'\la\phi,X_s\ra
       \prod_{i=1}^nH_i(X_{t_i})\big{\}}ds\cr
  &=&\lim_{k\rightarrow\infty}\mathbf{E}\big{\{}\la\phi,X_{t_{n+1}}^{(k)}\ra
     \prod_{i=1}^nH_i(X_{t_i}^{(k)})\big{\}}-
     \lim_{k\rightarrow\infty}\mathbf{E}\big{\{}\la\phi,X_{t_n}^{(k)}\ra
     \prod_{i=1}^nH_i(X_{t_i}^{(k)})\big{\}}\cr
  & &~~   -\lim_{k\rightarrow\infty}\int_{t_n}^{t_{n+1}}
      \mathbf{E}\big{\{}\mathcal{L}_k\la\phi,X_s^{(k)}\ra
       \prod_{i=1}^nH_i(X_{t_i}^{(k)})\big{\}}ds\cr
  &=&\lim_{k\rightarrow\infty}\mathbf{E}\big{\{}[\la\phi,X_{t_{n+1}}
  ^{(k)}\ra
  -\la\phi,X_{t_n}^{(k)}\ra-\int_{t_n}^{t_{n+1}}
    \mathcal{L}_k\la\phi,X_s^{(k)}\ra
    ds]\prod_{i=1}^nH_i(X_{t_i}^{(k)})\big{\}}\cr
  &=&0.
  \end{eqnarray*}
  Since $\{X_t:t\geq0\}$ is right continuous, the equality
  $$\mathbf{E}\big{\{}[\la\phi,X_{t_{n+1}}\ra-\la\phi,X_{t_n}\ra
  -\int_{t_n}^{t_{n+1}}
    \hat{\mathcal{L}}'\la\phi,X_s\ra ds]\prod_{i=1}^nH_i(X_{t_i})
    \big{\}}=0$$
  holds without the restriction $\{t_i\}_{i=1}^{n+1}\subset K(X)$.
  That is ($\ref{f3.15}$) is a martingale. Observe that
   if  $F(\mu)=f(\la1,\mu\ra)$ with $f\in C_0^2(\mbb R)$, then
   ${\cal A}_{\theta}F(\mu)=0$ and  $\cal
B_{\theta}F(\mu)$ is equal to
 \bgeqn
 \label{f3.17}
 \frac{\lz[\theta\wz\la1,\mu\ra]}{2\theta\la1,\mu\ra}
 \sum_{j=1}^{+\infty}
 (j-1)^2\la
 p_jf''(\la1,\mu)+\xi_j),\mu\ra
 \edeqn
 for some constant $0<\xi_j<(j-1)/\theta$. This follows from
 Taylor's expansion.  Recall that the sequence
$\{\lz_k\sigma_p^k/k\}$ and $\{\la1,\mu_k\ra\}$ are bounded.
 By the same argument as in the proof of Lemma $\ref{L2.1}$,
 we have
 $$
 \sup_k\mathbf{E}\la1,X_s^{(k)}\ra^2<\infty.
 $$
It follows from the Fatou's Lemma that ${\bf E}\la1,X_t\ra^2$ ia a
locally bounded function of $t\geq0$. Thus ($\ref{f3.15}$) is a
square-integrable martingale.

 \medskip

   \emph{Step 2}. We shall
show that under $\mathbf{Q}_{\mu}$
 \bgeqn
 \label{f3.18}
 \exp\{-\la\phi, \omega_t\ra\}-\exp\{-\la\phi,\omega_0\ra\}
 -\int_0^t\hat{\mathcal{L}}'\exp\{\la\phi,\omega_s\ra\}ds,~~
 t\geq0,
 \edeqn
is a martingale for $\phi\in C_{\partial}^2(\mathbb{R})^{++}$. Let
$\mu_k\in M_k(\mathbb{R})$ is given by
$$\mu_k=\frac{1}{k}\sum_{i=1}^{k\la1,\mu_k\ra}\delta_{x_i}.$$ Note that
 \beqlb
 \label{f3.19}
 \mathcal{A}_k\exp\{-\la\phi,\mu_k\ra\}
 \ar=\ar-\frac{1}{2}\exp\{-\la\phi,\mu_k\ra\}\la
 a\phi'',\mu_k\ra+\frac{1}{2k}\exp\{-\la\phi,\mu_k\ra\}\la
 (c\phi')^2,\mu_k\ra\cr
 \ar\ar+\frac{1}{2}\exp\{-\la\phi,\mu_k\ra\}
 \int_{\mathbb{R}^2}\rho(x-y)\phi'(x)\phi'(y)\mu_k(dx)\mu_k(dy)
 \eeqlb
and
 \beqlb
 \label{f3.20}
 \ar \ar\mathcal{B}_k\exp\{-\la\phi,\mu_k\ra\}\cr
 \ar=\ar\frac{k\lz_k(k\wedge\mu_k(1))}{\mu_k(1)}\bigg{\la}[
    \sum_{j=0}^{\infty}p_j(x)
    e^{-\la\phi,\mu_k\ra-\frac{j-1}{k}\phi(x)}-
    \sum_{j=0}^{\infty}p_j(x)
    e^{-\la\phi,\mu_k\ra}],\mu_k\bigg{\ra}\cr
 \ar=\ar\exp\{-\la\phi,\mu_k\ra\}\bigg{\la} \frac{k\lz_k(k\wedge\mu_k(1))}{\mu_k(1)}
    [\sum_{j=0}^{\infty}p_j(x)(e^{-\frac{j-1}{k}\phi(x)}-1)],\mu_k\bigg{\ra}\cr
 \ar=\ar\exp\{-\la\phi,\mu_k\ra\}\bigg{\la}\frac{(k\wedge\mu_k(1))}{\mu_k(1)}
 \psi_k(x,k-ke^{-\phi(x)/k})e^{\phi(x)/k},\mu_k\bigg{\ra}.
 \eeqlb
  Since for each
$h\geq0$, $\psi_k(x,z)\rightarrow\Psi_0(x,z)$ uniformly on
$\mathbb{R}\times[0,h]$,  we conclude  for $t\in K(X)$ a.s.
$\lim_{k\rightarrow\infty}\mathcal{L}_k\exp\{-\la
\phi,X_t^{(k)}\ra\}=\hat{\mathcal{L}}'\exp\{-\la \phi, X_t\ra\}$
boundedly by ($\ref{f3.16}$), ($\ref{f3.19}$), ($\ref{f3.20}$) and
the definition of $\hat{\mathcal{L}}'$. By the same argument as in
\textit{Step 1} we can get that ($\ref{f3.18}$) is a martingale.
That is
 \bgeqn
 W_t(\phi):=e^{-\la\phi,X_t\ra}-\int_0^t e^{-\la\phi,X_s\ra}[-
 \frac{1}{2}\la
 a\phi'',X_s\ra+\frac{1}{2}\int_{\hat{\mathbb{R}}}\la h(z-\cdot)
 \phi',X_s\ra^2dz+
 \la \Psi_0(\phi), X_s\ra]ds, ~t\geq0,
 \edeqn
 is a martingale with $\phi\in C^2_{\partial}(\mathbb{R})^{++}$, where $\Psi_0(\phi):=\Psi_0(x,\phi(x))$.
Then $\{\exp\{-\la\phi,X_t\ra\}:t\geq0\}$ is a special
semi-martingale with $\phi\in C^2_{\partial}(\mathbb{R})^{++}$.

 \emph{Step 3}.
 Let $S(\hat{\mathbb{R}})$ denote the space of finite signed Borel
 measures on $\hat{\mathbb{R}}$ endowed with the $\sigma$-algebra
 generated by the mappings $\mu\mapsto\la1,\mu\ra$ for all $f\in
 C(\hat{\mathbb{R}})$. Let
 $S(\hat{\mathbb{R}})^{\circ}=S(\hat{\mathbb{R}})\setminus\{0\}$.
 We define the optional random measure
 $N(ds,d\nu)$ on $[0,\infty)\times S(\hat{\mathbb{R}})^{\circ}$ by
 $$
 N(ds,d\nu)=\sum_{s>0}1_{\{\Delta
 X_s\neq0\}}\delta_{(s,\Delta X_s)}(ds,d\nu),
 $$
 where $\Delta X_s=X_s-X_{s-}\in S(\hat{\mathbb{R}})$. Let
 $\hat{N}(ds,d\nu)$ denote the predictable compensator of $N(ds,d\nu)$
 and let $\tilde{N}(ds,d\nu)$ denote the corresponding measure.
By the same argument as in the proof of Theorem $\ref{T2.3}$, we can
obtain that for $\phi\in C^2_{\partial}(\mbb R)$
\begin{eqnarray}
\label{f3.22}
 \la \phi,X_t\ra\ar=\ar\la \phi,\mu\ra+\int_0^t\la
 a\phi'',X_s\ra ds+M_t^c(\phi)+\int_0^{t+}\int_{S(\hat{\mathbb{R}})}
 \nu(\phi)\tilde{N}(ds,d\nu)\cr
 \ar\ar-
 \int_0^tds\int_{\hat {\mbb R}}X_s(dx)\phi(x)\int_l^{\infty}\xi
 \gamma(x,d\xi),
\end{eqnarray}
where $M_t^c(\phi)$ is a continuous local martingale. We also
conclude that the jump measure of the process $X$ has compensator
 \bgeqn
 \label{f3.23}
 \hat{N}(ds,d\nu)=dsX_s(dx)1_{\{0<\xi<l\}}\gamma(x,d\xi)\cdot\
 \dz_{\xi\dz_x}(d\nu),~~\nu\in
 M(\hat{\mathbb{R}})\setminus\{0\},
 \edeqn
and    for $\{\phi_i\}_{i=1}^2\subset
C_{\partial}^2(\mathbb{R})^{++}$,
 \begin{eqnarray}
 \label{f3.24}
 \la M^c(\phi_1),M^c(\phi_2)\ra_t
  &=&\frac{1}{2}\int_0^t
  \int_{\hat{\mathbb{R}}^2}\rho(x-y)\phi_1'(x)\phi_2'(y)
     X_s(dx)X_s(dy)ds\cr
  &&+\frac{1}{2}\int_0^t
  \int_{\hat{\mathbb{R}}^2}\rho(x-y)\phi_2'(x)\phi_1'(y)
     X_s(dx)X_s(dy)ds\cr
  &&+\int_0^t\la\sigma\phi_1\phi_2,X_s\ra ds.
  \end{eqnarray}

 \noindent Let $f\in C^2_0(\mathbb{R}^n)$ and
 $\{\phi_i\}_{i=1}^n\subset C^2_{\partial}(\mathbb{R})^{++}$.
  By ($\ref{f3.22}$),  ($\ref{f3.23}$), ($\ref{f3.24}$) and
  It\^{o}'s
  formula, we obtain
 \begin{eqnarray*}
 & &f(\la\phi_1,X_t\ra,\cdots,\la\phi_n,X_t\ra)\cr
 &=&f(\la\phi_1,X_0\ra,\cdots,\la\phi_n,X_0\ra)
    +\frac{1}{2}\sum_{i=1}^{n}
    \int_0^tf^i(\la\phi_1,X_s\ra,\cdots,\la\phi_n,X_s\ra)\la
    a\phi_i'',X_s\ra ds\cr
 \ar\ar+\frac{1}{2}\sum_{i,j=1}^n\int_0^t
 f^{ij}(\la\phi_1,X_s\ra,\cdots,\la\phi_n,X_s\ra)d\la
 M^c(\phi_i),M^c(\phi_j)\ra_t\cr
 \ar\ar
-\frac{1}{2}\sum_{i=1}^{n}\int_0^tds\int_{\hat{\mbb
R}}X_s(dx)\int_l^{\infty}\gamma(x,d\xi)
f^i(\la\phi_1,X_s\ra,\cdots,\la\phi_n,X_s\ra)\xi\phi_i(x)\cr
 \ar\ar+\int_0^tds\int_{\hat{\mathbb{R}}}X_s(dx)
 \int_0^{l}\gamma(x,d\xi)
 \{f(\la\phi_1,X_s\ra+\xi\phi_1(x),\cdots,\la\phi_n,X_s\ra
 +\xi\phi_n(x))\cr
 &  &~~~~-f(\la\phi_1,X_s\ra,\cdots,\la\phi_n,X_s\ra)
 -\xi\sum_{i=1}^n\phi_i(x)
 f^i(\la\phi_1,X_s\ra,\cdots,\la\phi_n,X_s\ra)\}\cr
 \ar\ar+(loc.mart.).
 \end{eqnarray*}
Hence
 $$F(X_t)-F(X_0)-\int_0^t\hat{\mathcal{L}}'F(X_s)ds,~~t\geq0,$$
 is a local martingale for each $F\in {\cal D}_0(\hat{\cal L}')$.
Since $f\in C_0^2(\mbb R^n)$ and $\phi_i\in C_{\partial}^2(\mbb
R)^{++}$, both $F$ and $\hat{\mathcal{L}}'F$ are bounded functions
on $M(\hat{\mbb {R}})$. Thus
 (\ref{f3.14}) is martingale. We complete the proof.\hfill$\Box$

\begin{lemma}
 \label{L3.3}
 Let $\mathcal{D}_0(\hat{\mathcal{L}}')$ be as in Lemma
 $\ref{L3.2}$. Then for each $\mu\in M(\hat{\mathbb{R}})$,
 there is a probability measure $\mathbf{Q}_{\mu}$ on
 $D([0,\infty),M(\hat{\mathbb{R}}))$ under which ($\ref{f3.14}$)
 is a martingale for each $F\in\mathcal{D}_0(\hat{\mathcal{L}}')$.
 \end{lemma}
 \textbf{Proof}. We only need to construct a sequence $\psi_k(x,z)$
 such that  for each
$h\geq0$, $\psi_k(x,z)\rightarrow\Psi_0(x,z)$ uniformly on
$\mathbb{R}\times[0,h]$, and $\frac{d}{dz}\psi_k(x,0+)\rightarrow
\frac{d}{dz}\Psi_0(x,0)$ uniformly on $\mbb R$ as
$k\rightarrow+\infty$. Moreover, $\{\frac{d^2}{dz^2}\psi_k(x,0+)\}$
should be a bounded sequence.

\noindent Let $\Psi_1(x,z)=\frac{1}{2}\sigma(x)z^2
  +\int_0^l(e^{-z\xi}-1+z\xi)\gamma(x,d\xi)$. We first define
   the sequences
$$\lz_{1,k}=1+k||\sigma||+
\sup_x\int_0^{l}\xi(1-e^{-k\xi})\gamma(x,d\xi)$$ and
$$g_{1,k}(x,z)=z+\frac{\Psi_1(x,k(1-z))}{k\lz_{1,k}}.$$
It is easy to check that $g_{1,k}(x,1)=1$ and
$$\frac{d^n}{dz^n}g_{1,k}(x,z)\geq0,~~~x\in\mathbb{R},~0\leq z\leq1,$$
for all integer $n\geq0$. Consequently, $g_{1,k}(x,\cdot)$ is a
probability generating function. Let $\psi_{1,k}(x,z)$ be defined by
($\ref{f3.13}$) with $(\lz_k, g_k)$ replaced by
$(\lz_{1,k},g_{1,k})$. Then
$$\psi_{1,k}(x,z)=\Psi_1(x,z)~~\textrm{for}~~0\leq z\leq k.$$
Let $b(x):=\int_l^{\infty}\xi\gamma(x,d\xi)$. Suppose $||b||>0$. Set
 $$g_{2,k}(x,z)=z+||b||^{-1}b(x)(1-z).
 $$
 Then $g_{2,k}(x,\cdot)$ is a probability generating function. Let
 $\lz_{2,k}=||b||$ and let $\psi_{2,k}(x,z)$ be defined by
 ($\ref{f3.13}$) with $(\lz_k,g_k)$ replaced by
 $(\lz_{2,k},g_{2,k})$. Then we have
 $$
 \psi_{2,k}(x,z)=b(x)z.
 $$
Finally we let $\lz_k=\lz_{1,k}+\lz_{2,k}$ and
$g_k=\lz_k^{-1}(\lz_{1,k}g_{1,k}+\lz_{2,k}g_{2,k})$. Then the
sequence $\psi_k$ defined by ($\ref{f3.13}$)  is equal to
$\psi_{1,k}+\psi_{2,k}$ which satisfies the required conditions
obviously.
 \qed
\smallskip
\begin{theorem}
\label{T3.1} Let $\{\omega_t:t\geq0\}$ denote the coordinate process
of $D([0,\infty),M(\mathbb{R}))$. Then for each $\mu\in
M(\mathbb{R})$ there is a probability measure $\mathbf{Q}_{\mu}$ on
$D([0,\infty),M(\mathbb{R}))$ such that $\{\omega_t:t\geq0\}$ under
$\mathbf{Q}_{\mu}$ is a solution of the
$(\mathcal{L}',\mu)$-martingale problem.
\end{theorem}
$\bf Proof.$  For each $\mu\in M(\mathbb{R})$, let
$\mathbf{Q}_{\mu}$ be the probability measure on
$D([0,\infty),M(\hat{\mathbb{R}}))$ provided by Lemma $\ref{L3.2}$.
We claim that  for any $T>0$
$$\mathbf{Q}_{\mu}\{\omega_t(\{\partial\})=0~\textrm{for all}~
 t\in[0,T]\}=1.$$ Consequently, $\mathbf{Q}_{\mu}$ is supported by
 $D([0,\infty),M(\mathbb{R}))$.
In fact, for any $\phi\in C_{\partial}^2(\mathbb{R})^{+}$, by
\textit{Step} 1 in the proof of Lemma $\ref{L3.2}$,  \bgeqn
\label{f3.25}M_t(\phi):=\la\phi,\omega_t\ra-\la\phi,\mu\ra-
\frac{1}{2}\int_0^t\la
 a\phi'',\omega_s\ra ds+
 \int_0^tds\int_{\hat {\mbb R}}\oz_s(dx)\phi(x)
 \int_l^{\infty}\xi\gamma(x,d\xi),~~t\geq0,
 \edeqn
  is a c\`{a}dl\`{a}g square-integrable martingale with quadratic
variation process given by
\begin{eqnarray*}
\la M(\phi)\ra_t
 =\int_0^t\big{\la}(\sigma+\int_0^{l}\xi^2\gamma(\cdot,d\xi))
 \phi^2,\omega_s\big{\ra}
ds+
 \int_0^tds\int_{\hat{\mathbb{R}}}\la
 h(z-\cdot)\phi',\omega_s\ra^2
 dz.
\end{eqnarray*}
 For
$k\geq1$, let
$$
 \phi_k(x)=\begin{cases}
 \exp\{-\frac{1}{|x|^2-k^2}\},& \textrm{if}~|x|>k,\\
  0,&\textrm{if}~|x|\leq k.
  \end{cases}
  $$
One can check that $\{\phi_k\}\subset C_{\partial}^2(\mathbb{R})$
such that  $\lim_{|x|\rightarrow\infty}\phi_k(x)=1$,
$\lim_{|x|\rightarrow\infty}\phi_k(x)'=0$ and
$\phi_k(\cdot)\rightarrow1_{\{\partial\}}(\cdot)$ boundedly and
pointwise. $||\phi_k'||\rightarrow0$ and $||\phi_k''||\rightarrow0$
as $k\rightarrow\infty$. Let
$\sigma_0=\sigma+\int_0^{l}\xi^2\gamma(\cdot,d\xi)$. By Theorem
1.6.10 of \cite{[IW89]}, we have
\begin{eqnarray*}
&&\mathbf{Q}_{\mu}\{\sup_{0\leq t\leq
T}|M_t(\phi_k)-M_t(\phi_j)|^2\}\cr
&\leq&4\int_0^T\mathbf{Q}_{\mu}\{\la\sigma_0
(\phi_k-\phi_j)^2,\omega_s\ra\}
ds+4\int_0^Tds\int_{\hat{\mathbb{R}}}\mathbf{Q}_{\mu}\{\la
 h(z-\cdot)(\phi_k'-\phi_j'),\omega_s\ra^2\}dz.\end{eqnarray*}
 By dominated convergence theorem, $\mathbf{Q}_{\mu}\{\sup_{0\leq
 t\leq
T}|M_t(\phi_k)-M_t(\phi_j)\ra|^2\}\rightarrow0$ as
$k,j\rightarrow0$. Therefore, there exists
$M^{\partial}=(M^{\partial}_t)_{t\geq0}$ such that for every $t>0$,
$$\mathbf{Q}_{\mu}\{|M_t(\phi_k)-M_t^{\partial}|^2\}
\rightarrow0
$$
and
$$
\sup_{0\leq s\leq t}|M_s(\phi_k)-M_s^{\partial}| \rightarrow0~~~
\textrm{in probabilty}
$$
as $k\rightarrow\infty$.  We obtain $M^{\partial}$ has
c\`{a}dl\`{a}g path. By Lemma 2.1.2 of \cite{[IW89]}, $M^{\partial}$
is a square-integrable martingale with zero mean. It follows from
($\ref{f3.25}$) that
$$M_t^{\partial}:=\omega_t(\{\partial\})+\int_0^tds\oz_s(\{\partial\}
)\int_l^{\infty}\xi\gamma(\partial,d\xi)$$ is a c\`{a}dl\`{a}g
square-integrable martingale with zero mean. Thus
$\mathbf{Q}_{\mu}(\omega_t(\{\partial\}))=0$. Then the claim follows
from the right continuity of
$\big{\{}\omega_t(\{\partial\}):t\geq0\big{\}}$. We have
$$
F(\omega_t)-F(\omega_0)-\int_0^t{\mathcal{L}'}F(\omega_s)ds,~~~t\geq0,
$$
is martingale for $F\in \mathcal{D}_0(\hat{\mathcal{L}}')$. Thus by
Remark \ref{rem2.1}, it is a local martingale for $F\in {\cal
D}(\cal L)$. \qed

\noindent Combining
 Theorem $\ref{T2.2}$ and Theorem $\ref{T3.1}$ we get
 that the $(\cal L',\mu)$-martingale problem is well-posed.
 The following theorem will show that the existence of solutions to
 $(\cal L,\mu)$-martingale problem.

\begin{theorem}
\label{T3.2}  For each $\mu\in M(\mathbb{R})$ there is a probability
measure $\bf Q_{\mu}$ on $(\Omega, \cal F)$  such that
  $\mathbf{Q}_{\mu}$ is a solution of the
$(\mathcal{L},\mu)$-martingale problem.
\end{theorem}
$\bf Proof.$ Let
$\lz_n(\mu)=1_{\{\la1,\mu\ra<n\}}\int\mu(dx)\int_l^{\infty}\gamma(x,d\xi)$
and define a transition function on $M(\mbb
R)\times\mathfrak{B}(M(\mbb R))$ by
$$\Gamma(\mu,d\nu):=\begin{cases}
\dz_{\mu}(d\nu),&
\int\mu(dx)\int_l^{\infty}\gamma(x,d\xi)=0,\\
(\int\mu(dx)\int_l^{\infty}\gamma(x,d\xi))^{-1}
\mu(dx)\gamma(x,d\xi)\dz_{\mu+\xi\dz_x}(d\nu),&otherwise.\end{cases}$$
Define ${\cal B}_n$ on $B(M(R))$ by
$$
B_nF(\mu):=\lz_n(\mu)\int(F(\nu)-F(\mu))\Gamma(\mu,d\nu)
=1_{\{\la1,\mu\ra<n\}}\int\mu(dx)\int_l^{\infty}(F(\mu+\xi\dz_x)-F(\mu))
\gamma(x,d\xi).
$$
Since the $(\cal L',\mu)$-martingale problem is well-posed, there
exists a semigroup $(Q_t')_{t\geq0}$ on $B(M(R))$ with transition
function given by ($\ref{f2.19}$) and full generator denoted by
${\cal L}'_0$. We can follow from Problem 4.11.3 of \cite{[EK86]} to
conclude that there exists a Markov process  denoted by
$X^n=\{X^n_t:t\geq0\}$ whose transition semigroup has full generator
given by  ${\cal L}'_0+{\cal B}_n$. In the following we assume that
$X^n_0=\mu~ a.s.$. Thus $({\cal L}'_0+{\cal B}_n,\mu)$-martingale
problem is well-posed. Since ${\cal L}'+{\cal B}_n\subset{\cal
L}'_0+{\cal B}_n$, $X^n$ is also a solution of $({\cal L}'+{\cal
B}_n,\mu)$-martingale problem. Let $U_n:=\{\mu\in M(\mbb
R):\la1,\mu\ra<n\}$. According to Theorem 4.3.6 of \cite{[EK86]},
there is a modification of $X^n$ with sample path in
$D([0,\infty),M(\mbb R))$.   Set
$$
\tau^n:=\inf\{t\geq0: \la1,X^n_t\ra\geq n\textrm{ or
}\la1,X^n_{t-}\ra\geq n\}
$$
and   $\tilde{X}^n=X^n_{\cdot\wz\tau^n}$. Then $\tilde{X}^n$ is a
solution of the stopped martingale problem for $({\cal L},U_n)$ and
by Theorem 4.6.1 of \cite{[EK86]}, $\tilde{X}^n$ is the unique
solution of the stopped martingale problem for $({\cal L}'_0+{\cal
B}_n,\dz_\mu,U_n)$. Put
$$
\tau^n_k:=\inf\{t\geq0:\la1,\tilde{X}_t^n\ra\geq k\textrm{ or }
\la1,\tilde{X}^n_{t-}\ra\geq k\}.
$$
For $k<n$, $\tilde{X}^n_{\cdot\wz\tau_k^n}$ is a solution of the
stopped martingale problem for $({\cal L}'_0+{\cal
B}_k,\dz_\mu,U_k)$ and hence has the same distribution as
$\tilde{X}^k$. On the other hand, since $\tilde{X}^n$ is a solution
of the stopped martingale problem for $({\cal L}, U_n)$, it follows
from Lemma $\ref{L2.3}$ that
$$
\sup_n{\bf E}\sup_{0\leq s\leq t}\la1,\tilde{X}^n_s\ra<\infty.
$$
Thus for each $t>0$,$$\lim_{n\rightarrow\infty}{\bf P}\{\tau^n\leq
t\}=0.$$   For any $k,~m\geq 1$, let $Y^k,~Y^m$ be two
$D([0,\infty),M(\mbb R))$-valued random variables such that they
have same distributions with $\tilde{X}^k$ and $\tilde{X}^m$
respectively and $Y^k(t)=Y^m(t)$ for $t\leq \tau^{k\wz m}$. Thus the
Skorohod distance between $Y^k$ and $Y^m$ is less than
$e^{-\tau^{k\wz m}}$. By Corollary 3.1.6 of \cite{[EK86]}, we
conclude that there exist a process $X^{\infty}$ such that
$\tilde{X}^n\Rightarrow X^{\infty}$. Let
 $$
 \tau^{\infty}_n=\inf\{t\geq0:\la1,X^{\infty}_t\ra\geq n\textrm{ or
}\la1,X^{\infty}_{t-}\ra\geq n\}.
 $$
Since the distribution of $\tilde{X}^m_{\cdot\wz\tau_n^m}$ does  not
depend on $m\geq n$, ${X}^{\infty}_{\cdot\wz\tau_n^{\infty}}$ has
the same distribution with $\tilde{X}^n$. Therefore, $${\bf
P}\{\tau^{\infty}_n\leq t\}={\bf P}\{\tau^n\leq t\}$$ and for each
$F\in \cal D(\cal L)$
$$
F({X}^{\infty}_{t\wz\tau^{\infty}_n})-\int_0^{t\wz\tau^{\infty}_n}{\cal
L}F({X}^{\infty}_s)ds
$$
is a martingale for each $n$. We see  $X^{\infty}$ is a solution of
the $({\cal L},\mu)$-martingale problem.\qed

\noindent Combining Theorem \ref{T2.6} and Theorem \ref{T3.2}, we
have that the $({\cal L},\mu)$-martingale problem is well-posed.
Thus we complete the construction of SDSM with general branching
mechanism.The next theorem gives another martingale characterization
of SDSM which is a direct consequence of Theorem \ref{T2.3} and
It\^{o}'s formula.

\begin{theorem} \label{T3.3} Let $\{\omega_t:t\geq0\}$ denote
the coordinate process of $D([0,\infty),M(\mathbb{R}))$. Then a
probability measure $\textbf{Q}_{\mu}$  on
$D([0,\infty),M(\mathbb{R}))$ is a solution of
$(\mathcal{L},\mu)$-martingale problem if and only if for $\mu\in
M(\mathbb{R})$ and $\phi\in C_c^2(\mathbb{R})^+$,
$\{\la\phi,\oz_t\ra\}$ is a semimartingale which has canonical
decomposition given by
 \begin{eqnarray}
\label{f3.26}
 \la \phi,\oz_t\ra\ar=\ar\la \phi,\mu\ra+\int_0^t\la
 a\phi'',\oz_s\ra ds+M_t^c(\phi)+\int_0^{t+}\int_{M(\mbb
R)^{\circ}}\la\phi,\nu\ra1_{\{\la1,\nu\ra<1\}}\tilde{N}(ds,d\nu)\cr
\ar\ar+\int_0^{t+}\int_{M(\mbb
R)^{\circ}}\la\phi,\nu\ra1_{\{\la1,\nu\ra \geq
1\}}N(ds,d\nu)-\int_0^{t}ds\int_{\mbb
R}\oz_s(dx)\int_1^{\infty}\xi\gamma(x,d\xi)\phi(x),
\end{eqnarray}
 where $\{M_t^c(\phi):t\geq0\}$ is a continuous local martingale with
 quadratic variation process given by
 \bgeqn
 \label{f3.27}
 \la M^c(\phi)\ra_t=\int_0^t\la\sigma\phi^2,\omega_s\ra ds+
 \int_0^tds\int_{\mathbb{R}}\la h(z-\cdot)\phi',\omega_s\ra^2
 dz,
 \edeqn
 and
 $$
 N(ds,d\nu)=\sum_{s>0}1_{\{\Delta
 \omega_s\neq0\}}\delta_{(s,\Delta \omega_s)}(ds,d\nu)
 $$
 is an optional random measure
  on $[0,\infty)\times M({\mathbb{R}})^{\circ}$,
 where $\Delta \oz_s=\oz_s-\oz_{s-}\in M({\mathbb{R}})$ and
 $\tilde{N}(ds,d\nu)$ denotes the corresponding martingale measure.
The predictable compensator of $N(ds,d\nu)$
 is given by
 $\hat{N}(ds,d\nu)=dsK(\oz_s,d\nu)$, where $K(\mu,d\nu)$
is determined  by
$$
\int_{M({\mathbb{R}})^{\circ}}F(\nu)K(\mu,d\nu)=
 \int_{\mbb R}\mu(dx)\int_0^{\infty}F(\xi\dz_x)\gamma(x,d\xi)
$$
for $F\in B(M(\mbb R))$.
 \end{theorem}

\section{Moment formulas, mean and spatial covariance measures}
In this section, we construct a dual process for SDSM and
investigate some properties of  SDSM. In accordance with the
notation used in Subsection 2.2, we can define a function-valued
Markov process by
 \bgeqn
 \label{4.1}Y'_t=P_{t-\tau_k}^{M_{\tau_k}}\Gamma_k
 P_{\tau_k-\tau_{k-1}}^{M_{\tau_{k-1}}}\Gamma_{k-1}\cdots
 P_{\tau_2-\tau_1}^{M_{\tau_1}}\Gamma_1P_{\tau_1}^{M_0}Y_0, \textrm{\
 \ \ }\tau_k\leq t<\tau_{k+1},~~0\leq k\leq M_0-1.
 \edeqn
 Let $X=\{X_t:t\geq0\}$ be an SDSM which is the unique solution of the
 martingale problem for $\cal L$.
If  for $m\geq2$,
$\sup_x[\int_0^{\infty}\xi^m\gamma(x,d\xi)]<\infty,$ then by the
same argument as in the proof of Lemma \ref{L2.1} and martingale
inequality, we have that
$${\bf E}\sup_{0\leq s\leq t}\la1,\oz_s\ra^m<\infty.$$ Then it follows from
the same argument of Theorem \ref{T2.1} that \bgeqn
 \label{4.2}
 \mathbf{E}\left\langle f,X_t^m\right\rangle
 =\mathbf{E}_{m,f}^{\sigma,\gamma}\big{[}\left\langle
 Y'_t,\mu^{M_t}\right\rangle
 \exp\big{\{}\int_0^t(2^{M_s}+\frac{M_s(M_s-1)}{2}-M_s-1)ds\big{\}}\big{]}
 \edeqn
for any $t\geq0$ and $f\in B(\mathbb{R}^m)$.

\noindent Skoulakis and Adler \cite{[SA01]} computed moments as a
limit of moments for the particle picture; see Section 3 of
\cite{[SA01]}. Stimulated by \cite{[SA01]}, in this section, we
compute moments via the dual relationship ($\ref{4.2}$). In fact, by
the construction ($\ref{4.1}$) of $\{Y'_t:t\geq0\}$ we have
\begin{eqnarray}
\label{dualequality}
 && \mathbf{E}_{m,f}^{\sigma,\gamma}\bigg{[}\la
 Y'_t,\mu^{M_t}\ra
 \exp\{\int_0^t(2^{M_s}+\frac{M_s(M_s-1)}{2}-M_s-1)ds\}\bigg{]}\cr
 &=&\la P_t^mf,\mu^m\ra\cr
 &&+\frac{1}{2}\sum_{i,j=1,i\neq j}^m\int_0^t
 \mathbf{E}_{m-1,\Psi_{ij}P_u^mf}
 ^{\sigma,\gamma}
 \bigg{[}\la
 Y'_{t-u},\mu^{M_{t-u}}\ra
 \exp\{\int_0^{t-u}(2^{M_s}+\frac{M_s(M_s-1)}{2}-M_s-1)ds\}\bigg{]}du\cr
 &&+\sum_{a=2}^m\bigg{(}\sum_{\{a\}}^m\int_0^t
 \mathbf{E}^{\sigma,\gamma}
 _{m-k+1,\Phi_{i_1,\cdots,i_a}P_u^mf}
  \big{[}\la
 Y'_{t-u},\mu^{M_{t-u}}\ra\cr
&&~~~~~~~~~~~~~~~~\times\exp\{\int_0^{t-u}
(2^{M_s}+\frac{M_s(M_s-1)}{2}-M_s-1)ds\}\big{]}du\bigg{)},
\end{eqnarray}
where $\{a\}=\{1\leq i_1<i_2<\cdots<i_a\leq m\}$. We remark that if
$\inf_{x\in\mathbb{R}}|c(x)|\geq\epsilon>0$, the semigroup
$(P_t^m)_{t>0}$ is uniformly elliptic and has density $p^m_t(x,y)$
satisfying
$$p_t^m(x,y)\leq\textrm{const}\cdot g_{\varepsilon
t}^m(x,y),~~t>0,x,y\in\mathbb{R}^m,$$ where $g_t^m(x,y)$ denotes the
transition density of the $m$-dimensional standard Brownian motion
(see Theorem 0.5 of \cite{[Dy65]}). In the following we always
assume that  $\sup_x[\int\xi^2\gamma(x, d\xi)]<\infty$.
\begin{theorem}
\label{T4.1} Suppose that $(\Omega, X_t, \mathbf{Q}_{\mu})$ is a
realization of the SDSM with parameters $(a,\rho,\Psi)$ with
$\inf_x|c(x)|\geq\epsilon>0$. Let $f\in B(\mathbb{R})$ and $t>0$.
Then we have the first moment formula for $X$ as follows:
 \bgeqn
 \label{firstmo}
 \mathbf{E}(\la
 f,X_t\ra)=\int_{\mathbb{R}}\int_{\mathbb{R}}f(y)p_t(x,y)dy\mu(dx),
 \edeqn
and   $\forall\, 0<s\leq t$, $f\in B(\mathbb{R})$ and $g\in
B(\mathbb{R})$, we have the second order moment formula
 \begin{eqnarray}
 \label{secondm}
 &&\mathbf{E}(\la f,X_{s}\ra\la g,X_{t}\ra)\cr
 &=&\mathbf{E}(\la f, X_{s}\ra\la P_{t-s}g,X_{s}\ra)\cr
 &=&\int_{\mathbb{R}}\int_{\mathbb{R}}\int_{\mathbb{R}^2}f(y_1)
  \left(\int_{\mathbb{R}}g(z)p_{t-s}(y_2,z)dz\right)
   p_{s}^2(x,y;y_1,y_2)dy_1dy_2\mu(dy)\mu(dx)\cr
 &&+\int_0^{s}du\int_{\mathbb{R}}\mu(dx)\int_{\mathbb{R}}dy
  \int_{\mathbb{R}^2}dy_1dy_2p_{s-u}(x,y)\sigma(y)p_u^2(y,y;y_1,y_2)\cr
  &&~~\times
  f(y_1)\left(\int_{\mathbb{R}}p_{t-s}(y_2,z)g(z)dz\right)\cr
  &&+\int_0^{s}du\int_{\mathbb{R}}\mu(dx)\int_{\mathbb{R}}dy
  \int_{\mathbb{R}^2}dy_1dy_2p_{s-u}(x,y)
  \left(\int_0^{\infty}\xi^2\gamma(y,d\xi)\right)p_u^2(y,y;y_1,y_2)\cr
  &&~~\times
  f(y_1)\left(\int_{\mathbb{R}}p_{t-s}(y_2,z)g(z)dz\right).
 \end{eqnarray}\end{theorem}
\textbf{Proof}. ($\ref{firstmo}$) is a direct conclusion of
($\ref{dualequality}$). Using ($\ref{firstmo}$) and  Markov property
of $X$ we have $\mathbf{E}(\la f,X_{s}\ra\la g,X_{t}\ra)
 =\mathbf{E}(\la f, X_{s}\ra\la P_{t-s}g,X_{s}\ra)$. Then
($\ref{secondm}$) is also a direct conclusion of
($\ref{dualequality}$).\qed

\noindent Following \cite{[SA01]}, we define two deterministic
measures as follows:
\begin{enumerate}
\item The \emph{mean measure} $m_t$ defined on
$\mathcal{B}(\mathbb{R})$ by
$$
m_t(A)=\mathbf{E}(X_t(A)).
$$
\item The \emph{spatial measure} $s_t$ defined on
$\mathcal{B}(\mathbb{R}\times\mathbb{R})$ by
$$
s_t(A_1\times A_2)=\mathbf{E}(X_t(A_1)X_t(A_2)).
$$
\end{enumerate}
By Theorem $\ref{T4.1}$, we have following proposition.
\begin{proposition}
For all $t>0$ the measures $m_t$ and $s_t$ have densities with
respect to Lebesgue measure,denoted by $m(t;y)$ and $s(t;y_1,y_2)$,
respectively. We have that
$$
m(t;y)=\int_{\mathbb{R}}p_t(x,y)\mu(dx)
$$
for all $y\in\mathbb{R}$ and
 \begin{eqnarray}
 \label{density}
 s(t;y_1,y_2)
 &=&\int_{\mathbb{R}^2}p^2_t(y,z;y_1,y_2)\mu(dy)\mu(dz)\cr
 &&~~~+\int_0^tds\int_{\mathbb{R}}\mu(dy)
 \int_{\mathbb{R}}dz\sigma(z)p_s^2(z,z;y_1,y_2)p_{t-s}(y,z)\cr
 &&~~~+\int_0^tds\int_{\mathbb{R}}\mu(dy)
 \int_{\mathbb{R}}dz
 \int_0^{\infty}\xi^2\gamma(z,d\xi)p_s^2(z,z;y_1,y_2)p_{t-s}(y,z)
 \end{eqnarray}
for all $y_1,y_2\in\mathbb{R}$.
\end{proposition}

\bigskip
\textbf{Acknowledgement}. I would like to give my sincere thanks to
my supervisor Professor Zenghu Li for his encouragement and helpful
discussions. I also would like to express my sincere gratitude to
the referee and the AE for their encouragement and useful comments.
I thank Professor Mei Zhang for her useful suggestions to the
earlier version of the paper.

\end{document}